\documentclass[12pt]{article}
\usepackage{amssymb,amsmath}
\begin{document}
\font\Goth=eufm10 scaled 1440 \font\SY=msam10
\def\ccc{\hbox{{\Goth c}}}
\font\VV=cmsy10 scaled 1728
\def\vvv{\mathop{\hbox{\VV\char"5F}}}
\def\N{{\mathbb N}}
\def\Z{{\mathbb Z}}
\def\R{{\mathbb R}}
\def\C{{\mathbb C}}
\def\T{{\mathbb T}}
\def\K{{\mathbb K}}
\def\U{\overline{U}}
\def\D{{\mathbb D}}
\def\zp{\Z_+}
\def\dim{{\rm dim}\,}
\def\diam{{\rm diam}\,}
\def\square{\hbox{\SY\char"03}}
\def\epsilon{\varepsilon}
\def\phi{\varphi}
\def\kappa{\varkappa}
\def\theorem#1{\smallskip\noindent
{\scshape Theorem} {\bf #1}{\bf .}\hskip 8pt\sl}
\def\defin#1{\smallskip\noindent
{\scshape Definition} {\bf #1}{\bf .}\hskip 6pt}
\def\prop#1{\smallskip\noindent
{\scshape Proposition} {\bf #1}{\bf .}\hskip 8pt\sl}
\def\lemma#1{\smallskip\noindent
{\scshape Lemma} {\bf #1}{\bf .}\hskip 6pt\sl}
\def\cor#1{\smallskip\noindent
{\scshape Corollary} {\bf #1}{\bf .}\hskip 6pt\sl}
\def\quest#1{\smallskip\noindent
{\scshape Question} {\bf #1}{\bf .}\hskip 6pt\sl}
\def\epr{\smallskip\rm}
\def\rem#1{\smallskip\noindent
{\scshape Remark}  {\bf #1}.\hskip 6pt}
\def\ii{\mathrel{\Longrightarrow}}
\def\wz{\thinspace}
\def\proof{P\wz r\wz o\wz o\wz f.\hskip 6pt}
\def\leq{\leqslant}
\def\geq{\geqslant}
\def\RR{{\cal R}}
\def\ssub#1#2{#1_{{}_{{\scriptstyle #2}}}}
\def\li{\,\text{\rm span}\,}
\def\V{{\cal V}}
\def\H{{\cal H}}
\def\slim{\mathop{\hbox{$\overline{\hbox{\rm lim}}$}}}
\def\ilim{\mathop{\hbox{$\underline{\hbox{\rm lim}}$}}}

\title{On similarity of quasinilpotent operators}

\author{S.~Shkarin}

\date{}

\maketitle

\smallskip

\leftline{King's College London, Department of Mathematics}
\leftline{Strand, London WC2R 2LS, UK} \leftline{\bf e-mail: \tt
stanislav.shkarin@kcl.ac.uk}

\rm\normalsize

\begin{abstract}

Bounded linear operators on separable Banach spaces algebraically
similar to the classical Volterra operator $V$ acting on $C[0,1]$
are characterized. From this characterization it follows that $V$
does not determine the topology of $C[0,1]$, which answers a
question raised by Armando Villena. A sufficient condition for  an
injective bounded linear operator on a Banach space to determine its
topology is obtained. From this condition it follows, for instance,
that the Volterra operator acting on the Hardy space $\H^p$ of the
unit disk determines the topology of $\H^p$ for any
$p\in[1,\infty]$.



\end{abstract}

\section{Introduction}

All vector spaces in this article are assumed to be over the field
$\C$ of complex numbers unless specified otherwise. The following
problem falls into the family of questions on the so-called
automatic continuity. Suppose that $T$ is a bounded linear operator
acting on a Banach space $X$ and $\ssub{\|\cdot\|}{1}$ is another
complete norm on $X$, with respect to which $T$ is also bounded.
Does it follow that $\ssub{\|\cdot\|}{1}$ is equivalent to the
initial norm? According to the Banach Inverse Mapping Theorem
\cite{rud}, this question is equivalent to the following one. Should
$\ssub{\|\cdot\|}{1}$ be continuous? If the answer is affirmative,
we say that $T$ {\it determines the topology of} $X$. One may ask
the same question about a family of operators. There are several
results both positive and negative for different type of operators
and families of operators, see for instance,
\cite{alexvi,exmevi,ja1,ja2,vi3,vi1,vi2,vi4}. The following question
was raised by Armando Villena in 2000; it can be found in the
Belfast Functional Analysis Day problem book.

\quest1 Does the Volterra operator
$$
\smash{V:C[0,1]\to C[0,1],\qquad Vf(x)=\int\limits_0^x f(t)\,dt}
$$
determine the topology of $C[0,1]$?\epr

We answer this question negatively. To this end bounded linear
operators on Banach spaces (algebraically) similar to $V$ are
characterized. Two linear operators $T:X\to X$ and $S:Y\to Y$ are
said to be {\it similar} if there exists an invertible linear
operator $G:X\to Y$ such that $GT=SG$. Note that we consider
similarity {\bf only in the algebraic sense}, that is, even if $X$
and $Y$ carry some natural topologies, with respect to which $T$ and
$S$ are continuous, the operator $G$ {\bf is not assumed to be
continuous}.

\prop{1.1} Let $T$ be a bounded linear operator on a Banach space
$X$. Then the following conditions are equivalent.
\begin{itemize}
\item[{\bf (T1)}]$T$ determines the topology of $X$;
\item[{\bf (T2)}]if $S$ is a bounded linear operator acting on a
Banach space $Y$, $S$ is similar to $T$ and $G:X\to Y$ is an
invertible linear operator such that $GT=SG$, then $G$ is bounded.
\end{itemize}
\rm

\proof If $T$ does not determine the topology of $X$, then there
exists a complete norm $\ssub{\|\cdot\|}0$ on $X$ not equivalent to
the initial norm and such that $T$ is bounded with respect to
$\ssub{\|\cdot\|}0$. Let $G$ be the identity operator $Gx=x$ from
$X$ endowed with the initial norm $\ssub{\|\cdot\|}X$ to $X$ endowed
with the norm $\ssub{\|\cdot\|}0$. Since both norms are complete and
not equivalent, $G$ is unbounded. If $S$ is the operator $T$ acting
on the Banach space $(X,\ssub{\|\cdot\|}0)$, then the equality
$GT=SG$ is obviously satisfied. Thus, (T2) implies (T1).

Suppose now that (T2) is not satisfied. Pick a bounded linear
operator $S$ on a Banach space $Y$ and an unbounded invertible
linear operator $G:X\to Y$ such that $GT=SG$. Consider the norm on
$X$ defined by the formula $\ssub{\|x\|}0=\ssub{\|Gx\|}{Y}$. Then
$G$ is an isometry from $(X,\ssub{\|\cdot\|}{0})$ onto the Banach
space $Y$. Hence $X$ with the norm $\ssub{\|\cdot\|}{0}$ is
complete. Since $G$ is unbounded, the norm $\ssub{\|\cdot\|}{0}$ is
not equivalent to the initial one. On the other hand the operator
$T$ acting on $(X,\ssub{\|\cdot\|}{0})$ is bounded since it is
isometrically similar to $S$. Thus, $T$ does not determine the
topology of $X$. \square

We shall give a negative answer to Question~1 by means of
Proposition~1.1, proving that the Volterra operator acting on
$C[0,1]$ is similar to a bounded linear operator acting on a Banach
space non-isomorphic to $C[0,1]$. Let $T$ be a linear operator
acting on a linear space $X$. The {\it spectrum} of $T$ is
$$
\sigma(T)=\{z\in\C:T-zI\ \ \text{is non-invertible}\}.
$$
According to the Banach inverse mapping theorem, if $X$ is a
Fr\'echet space and $T$ is continuous, then $\sigma(T)$ coincides
with the conventional spectrum since $(T-zI)^{-1}$ is continuous
whenever $z\in\C\setminus\sigma(T)$. This is not true for continuous
linear operators on general topological vector spaces. We say that
$T$ is {\it quasinilpotent} if $\sigma(T)\subseteq\{0\}$. Everywhere
below $\Z$ is the set of integers and $\zp$ is the set of
non-negative integers.

The following theorem characterizes similarity to the Volterra
operator.

\theorem{1.2}Let $T$ be a bounded linear operator on a Banach space
$X$ of algebraic dimension $\ccc=2^{\aleph_0}$. Then the following
conditions are equivalent.
\begin{itemize}
\item[{\bf (C1)}]$T$ is similar to the Volterra operator $V$ acting
on $C[0,1]$;
\item[{\bf (C2)}]$T$ is injective, quasinilpotent and satisfies the
closed finite descent condition, that is, there exists $m\in\zp$ for
which
\begin{equation}
\overline{T^m(X)}=\overline{T^{m+1}(X)}. \label{t121}
\end{equation}
\end{itemize}
\rm

Clearly (\ref{t121}) is satisfied for $m=0$ if the range of $T$ is
dense. Since the algebraic dimension of any separable infinite
dimensional Banach space is $\ccc$, we have

\cor{1.3}Any injective quasinilpotent bounded linear operator $T$
with dense range acting on a separable Banach space is similar to
the Volterra operator acting on $C[0,1]$. \epr

\theorem{1.4}The Volterra operator does not determine the topology
of $C[0,1]$. \epr

\proof Let $V$ be the Volterra operator acting on $C[0,1]$ and $V_2$
be the same operator acting on $L_2[0,1]$. By Theorem~1.2 $V$ and
$V_2$ are similar. Since Banach spaces $C[0,1]$ and $L_2[0,1]$ are
not isomorphic, the similarity operator can not be bounded.
Proposition~1.1 implies now that $V$ does not determine the topology
of $C[0,1]$. \square

We also provide a new sufficient condition for an injective bounded
linear operator on a Banach space to determine its topology. It
follows that certain injective quasinilpotent operators do determine
the topology of the Banach space on which they act. In particular,
it is true for the Volterra operator acting on the Hardy space
$\H^p$ of the unit disk for any $p\in[1,\infty]$.

\theorem{1.5}Let $T$ be a bounded injective linear operator on a
Banach space $X$ such that
\begin{align}
&\bigcap_{n=0}^\infty \overline{T^n(X)}=\{0\}\ \ \text{and}
\label{151}
\\
&\text{there exists $n\in\zp$ such that $\overline{T^{n+1}(X)}$ has
finite codimension in $\overline{T^n(X)}$.} \label{152}
\end{align}
Then $T$ determines the topology of $X$. \epr

In Section~2 linear operators with empty spectrum are characterized
up to similarity. In Section~3 we introduce the class of tame
operators and characterize up to similarity tame injective
quasinilpotent operators. Section~4 is devoted to auxiliary lemmas
with the help of which tame injective bounded operators on Banach
spaces are characterized in Section~5. Theorem~1.2 and Theorem~1.5
are proved in Section~6 and Section~7 respectively. In Section~8 of
concluding remarks we discuss the previous results and raise few
problems.

\section{Operators with empty spectrum}

In what follows $\RR$ stands for the field of complex rational
functions considered also as a complex vector space, and
$$
M:\RR\to\RR,\quad Mf(z)=zf(z).
$$

\lemma{2.1}Let $T$ be a linear operator on a linear space $X$ such
that $\sigma(T)=\varnothing$. Consider the multiplication operation
from $\RR\times X$ to $X$ defined by the formula $r\cdot x=r(T)x$.
Then this multiplication extends the natural multiplication by
complex numbers and turns $X$ into a linear space over the field
$\RR$. \epr

\proof Since $\sigma(T)=\varnothing$, $r(T)$ is a well-defined
linear operator on $X$ for any $r\in\RR$ and $r(T)$ is invertible if
$r\neq 0$. Thus, the operation $(r,x)\mapsto r\cdot x$ is
well-defined. The verification of axioms of the vector space is
fairly elementary. \square

\theorem{2.2}Let $T$ be a linear operator acting on a linear space
$X$ such that $\sigma(T)=\varnothing$. Then there exists a unique
cardinal $\mu=\mu(T)$ such that $T$ is similar to the direct sum of
$\mu$ copies of the multiplication operator $M$. Moreover, $\mu(T)$
coincides with the algebraic dimension of $\ssub X\RR$, being $X$,
considered as a linear space over $\RR$ with the multiplication
$r\cdot x=r(T)x$. In particular, two linear operators $T:X\to X$ and
$S:Y\to Y$ with $\sigma(T)=\sigma(S)=\varnothing$ are similar if and
only if $\mu(T)=\mu(S)$. \epr

\proof Let $\nu$ be the algebraic dimension of $\ssub{X}{\RR}$ and
$\{x_\alpha:\alpha\in A\}$ be a Hamel basis in $\ssub{X}{\RR}$. Then
the cardinality of $A$ is $\nu$. For each $\alpha\in A$ let
$X_\alpha$ be the $\RR$-linear span of the one-element set
$\{x_\alpha\}$. Then $X_\alpha$ are linear subspaces of
$\ssub{X}{\RR}$ and  therefore they are $\C$-linear subspaces of
$X$. Moreover $X$ is the direct sum of $X_\alpha$. Since for any
$x\in X$, the vectors $x$ and $Tx$ are $\RR$-collinear, we see that
each $X_\alpha$ is $T$-invariant. Moreover, for any $x\in X$ and any
$r\in\RR$, $T(r\cdot x)=Mr\cdot x$ and therefore, for any $\alpha\in
A$, the restriction $T_\alpha$ of $T$ to $X_\alpha$ is similar to
$M$ with the similarity provided by the operator $G_\alpha:\RR\to
X_\alpha$, $G_\alpha r=r\cdot x_\alpha$. Thus, $T$ is similar to the
direct sum of $\nu$ copies of $M$.

Suppose now that $\mu$ is a cardinal and $T$ is similar to the
direct sum of $\mu$ copies of $M$. Then there exists a set $A$ of
cardinality $\mu$ and a family $\{X_\alpha:\alpha\in A\}$ of
$T$-invariant linear subspaces of $X$ such that $X$ is the direct
sum of $X_\alpha$ and for each $\alpha\in A$, the restriction
$T_\alpha$ of $T$ to $X_\alpha$ is similar to $M$. Let
$G_\alpha:\RR\to X_\alpha$ be a $\C$-linear invertible linear
operator such that $G_\alpha M=T_\alpha G_\alpha$ and
$x_\alpha=G_\alpha(1)$. One can easily verify that
$G_\alpha(r\rho)=r\cdot G_\alpha(\rho)$ for any $r,\rho\in \RR$.
Hence for any pairwise different $\alpha_1,\dots,\alpha_n\in A$, the
$\RR$-linear span of $\{x_{\alpha_1},\dots,x_{\alpha_n}\}$ coincides
with the direct sum of $X_{\alpha_1},\dots,X_{\alpha_n}$. Therefore
$\{x_\alpha:\alpha\in A\}$ is a Hamel basis in $\ssub{X}{\RR}$.
Thus, $\mu=\nu$. \square

\section{Similarity of tame operators}

For a linear operator $T$ on a linear space $X$ we denote
$$
\ssub XT=\bigcap_{n=0}^\infty T^n(X).
$$
Clearly $\ssub XT$ is a $T$-invariant linear subspace of $X$. Let
$T_0:\ssub XT\to \ssub XT$ be the restriction of $T$ to the
invariant subspace $\ssub XT$.

\lemma{3.1}Let $T$ be an injective quasinilpotent operator on a
linear space $X$. Then $\sigma(T_0)=\varnothing$.\epr

\proof Let $z\in \C$. We have to prove that $T_0-zI$ is invertible.
Since $T_0-zI$ is the restriction of the injective operator $T-zI$,
we see that $T_0-zI$ is injective. It remains to verify surjectivity
of $T_0-zI$. Let $x_0\in \ssub XT$. Since $T$ is injective and
$x_0\in \ssub XT$, for any $n\in\zp$ there exists a unique $x_n\in
X$ such that $T^nx_n=x_0$. From injectivity of $T$ it also follows
that $x_n=T^mx_{m+n}$ for each $m,n\in\zp$. Therefore all $x_n$
belong to $\ssub XT$.

{\bf Case } $z=0$. Since $x_1\in X_0$ and $Tx_1=T_0x_1=x_0$, we see
that $T_0=T_0-zI$ is surjective.

{\bf Case } $z\neq 0$. Then $T-zI$ is invertible. Denote
$y=(T-zI)^{-1}x_0$. Then for any $n\in\zp$,
$$
y=(T-zI)^{-1}T^nx_n=T^nw_n, \ \ \text{where}\ \ w_n=(T-zI)^{-1}x_n.
$$
Hence $y\in \ssub XT$ and $(T_0-zI)y=x_0$. Thus, $T_0-zI$ is
surjective. \square

In order to formulate the main result of this section we need some
additional notation. If $X$ is a linear space, $E$ is a linear
subspace of $X$ and $x,y\in X$, we write
$$
x\equiv y\ (\bmod\, E)
$$
if $x-y\in E$. We also say that a family $\{x_\alpha\}_{\alpha\in
A}$ of elements of $X$ is {\it linearly independent modulo} $E$ if
for any pairwise different $\alpha_1,\dots,\alpha_n\in A$ and any
complex numbers $c_1,\dots,c_n$, the inclusion
$c_1x_{\alpha_1}+{\dots}+c_nx_{\alpha_n}\in E$ implies $c_j=0$ for
$1\leq j\leq n$. Clearly linear independence of
$\{x_\alpha\}_{\alpha\in A}$ modulo $E$ is equivalent to linear
independence of $\{\pi(x_\alpha)\}_{\alpha\in A}$ in $X/E$, where
$\pi:X\to X/E$ is the canonical map.

\defin1 We say that a linear operator $T$ acting on a linear space
$X$ is {\it tame} if for any sequence $\{x_n\}_{n\in\zp}$ of
elements of $X$ there exists $x\in X$ such that for each $n\in\zp$,
\begin{equation}
x\equiv \sum_{k=0}^{n}T^kx_k\ (\bmod\, T^{n+1}(X)). \label{tame1}
\end{equation}

Obviously surjective operators are tame and tameness is invariant
under similarities. It is not clear {\it a priori} whether there are
injective tame operators $T$ with $\sigma(T)=\{0\}$. We shall show
in the next section that the Volterra operator acting on $C[0,1]$ is
tame, which motivates the study of similarity of injective
quasinilpotent tame operators. In what follows $\ssub\dim\K X$
stands for the algebraic dimension of the linear space $X$ over the
field $\K$.

\theorem{3.2}Let $T:X\to X$ and $S:Y\to Y$ be two injective
quasinilpotent tame operators. Then $T$ and $S$ are similar if and
only if $\ssub\dim\C X/T(X)=\ssub\dim\C Y/T(Y)$ and
$\ssub{\dim}{\RR}\ssub XT=\ssub{\dim}{\RR}\ssub YS$, where the
multiplication of $x\in\ssub XT$ and $y\in\ssub YS$ by $r\in\RR$ are
given by $r\cdot x=r(T_0)x$ and $r\cdot y=r(S_0)y$ respectively.
\epr

The rest of the section is devoted to the proof of Theorem~3.2. We
need some preparation.

\subsection{Auxiliary Lemmas}

We set
$$
\RR_+=\{f\in\RR: f(0)\neq\infty\}.
$$
If $T$ is injective and quasinilpotent, then $r(T)$ is a
well-defined injective linear operator on $X$ for any
$r\in\RR_+\setminus\{0\}$ and $r(T)$ is invertible if
$r(0)\in\C\setminus\{0\}$.

\defin2 Let $T$ be an injective quasinilpotent operator acting on
a linear space $X$. We say that a family of vectors
$\{x_\alpha\}_{\alpha\in A}$ in $X$ is $T$-{\it independent} if for
any pairwise different $\alpha_1,\dots,\alpha_n\in A$ and any
$r_1,\dots,r_n\in\RR_+$, the inclusion
$r_1(T)x_{\alpha_1}+{\dots}+r_n(T)x_{\alpha_n}\in\ssub XT$ implies
$r_j=0$ for $1\leq j\leq n$.

From Lemma~3.1 it follows that if $T:X\to X$ is injective and
quasinilpotent, then $r(T)(\ssub XT)=\ssub XT$ for any
$r\in\RR_+\setminus\{0\}$. From this observation it follows that a
family $\{x_\alpha\}_{\alpha\in A}$ in $X$ is $T$-{\it independent}
if and only if for any pairwise different
$\alpha_1,\dots,\alpha_n\in A$ and any polynomials $p_1,\dots,p_n$
the inclusion $p_1(T)x_{\alpha_1}+{\dots}+p_n(T)x_{\alpha_n}\in\ssub
XT$ implies $p_j=0$ for $1\leq j\leq n$. Applying the Zorn Lemma to
the set of $T$-independent families partially ordered by inclusion,
we see that there are maximal $T$-independent families.

\lemma{3.3}Let $T$ be an injective quasinilpotent operator acting on
a linear space $X$ and $\{x_\alpha\}_{\alpha\in A}$ be a maximal
$T$-{\it independent} family. Then for any $x\in X$ there exists a
unique  finite $($maybe empty$)$ subset $\Lambda=\Lambda(x)$ of $A$
such that
\begin{equation}
T^nx=w_n+\sum_{\alpha\in\Lambda}r_{\alpha,n}(T)x_\alpha \label{tin}
\end{equation}
for some $n\in\zp$, $w_n\in \ssub XT$ and non-zero
$r_{\alpha,n}\in\RR_+$. Moreover $w_n$ and $r_{\alpha,n}$ are
uniquely determined by $n$ and $x$. Finally, if the decomposition
$(\ref{tin})$ does exist for some $n$ then it exists for all greater
$n$ and $w_m=T^{m-n}w_n$, $r_{\alpha,m}=M^{m-n}r_{\alpha,n}$ for
$m\geq n$. \epr

\proof Let $B$ and $C$ be finite subsets of $A$, $r_\alpha$ for
$\alpha\in B$ and $\rho_\alpha$ for $\alpha\in C$ be non-zero
elements of $\RR_+$, $m,n\in\zp$ and $u,v\in\ssub XT$ be such that
$$
T^nx=u+\sum_{\alpha\in B}r_\alpha(T)x_\alpha\ \ \text{and}\ \
T^mx=v+\sum_{\alpha\in C}\rho_\alpha(T)x_\alpha.
$$
Let also $D=B\cup C$ and $r'_\alpha$, $\rho'_\alpha\in\RR_+$ for
$\alpha\in D$ be defined by the formulas: $r'_\alpha=M^mr_\alpha$
for $\alpha\in B$, $\rho'_\alpha=M^n\rho_\alpha$ for $\alpha\in C$,
$r'_\alpha=0$ for $\alpha\in C\setminus B$ and $\rho'_\alpha=0$ for
$\alpha\in B\setminus C$. Then
\begin{equation}
T^{m+n}x=T^mu+\sum_{\alpha\in D}r'_\alpha(T)x_\alpha=
T^nv+\sum_{\alpha\in D}\rho'_\alpha(T)x_\alpha. \label{tmn}
\end{equation}
Hence
$$
\sum_{\alpha\in D}(r'_\alpha-\rho'_\alpha)(T)x_\alpha\in\ssub XT.
$$
Since $\{x_\alpha\}_{\alpha\in A}$ is $T$-{\it independent}, we have
$\rho'_\alpha=r'_\alpha$ for each $\alpha\in D$. Hence $B=C=D$ and
$M^mr_\alpha=M^n\rho_\alpha$ for any $\alpha\in D$. Substituting
these equalities into (\ref{tmn}), we see that $T^mu=T^nv$. This
means that the set $\Lambda$ is uniquely determined by $x$ and
$w_n$, $r_{\alpha,n}$ are uniquely determined by $x$ and $n$.
Moreover, if (\ref{tin}) is satisfied then it is satisfied for
greater $n$ with $w_m=T^{m-n}w_n$ and
$r_{\alpha,m}=M^{m-n}r_{\alpha,n}$ for $m\geq n$.

Suppose now that $x\in X$ does not admit a decomposition
(\ref{tin}). It easily follows that $\{x\}\cup
\{x_\alpha\}_{\alpha\in A}$ is $T$-independent, which contradicts
the maximality of $\{x_\alpha\}_{\alpha\in A}$. \square

\lemma{3.4}Let $T$ be an injective linear operator acting on a
linear space $X$ and $E$ be a linear subspace of $X$ such that
$E\oplus T(X)=X$. Then for any $x\in X$ there exists a unique
sequence $\{x_n\}_{n\in\zp}$ of elements of $E$ such that
$(\ref{tame1})$ is satisfied for any $n\in\zp$. \epr

\proof Let $x\in X$. It suffices to prove that for any $m\in\zp$,
there exist unique $x_0,\dots,x_m\in E$ for which $(\ref{tame1})$ is
satisfied for $n\leq m$. We achieve this using induction with
respect to $m$.

Since $E\oplus T(X)=X$, there exists a unique $x_0\in E$ for which
$x-x_0\in T(X)$. This inclusion is exactly (\ref{tame1}) for $n=0$.
The basis of induction is constructed. Let now $m$ be a positive
integer. Assume that there exist unique $x_0,\dots,x_{m-1}\in E$ for
which $(\ref{tame1})$ is satisfied for $n\leq m-1$. Since $T$ is
injective, $(\ref{tame1})$ for $n=m-1$ implies that there exists a
unique $u\in X$ such that $x-x_0-Tx_1-{\dots}-Tx_{m-1}=T^mu$. Since
$E\oplus T(X)=X$, there exists a unique $x_m\in E$ for which
$u-x_m\in T(X)$. Since the last inclusion is equivalent to
$(\ref{tame1})$ for $n=m$, the induction step is complete and so is
the proof of the lemma. \square

\subsection{Proof of Theorem~3.2}

Let $T$ and $S$ be similar. Since the algebraic codimension of the
range of a linear operator is a similarity invariant, we have
$\ssub\dim\C X/T(X)=\ssub\dim\C Y/T(Y)$. Let $G:X\to Y$ be an
invertible linear operator such that $GT=SG$. Clearly $G(\ssub
XT)=\ssub YS$. Therefore $T_0$ and $S_0$ are similar. According to
Lemma~3.1 and Theorem~2.2 $\ssub{\dim}{\RR}\ssub
XT=\ssub{\dim}{\RR}\ssub YS$.

Suppose now that $\ssub\dim\C X/T(X)=\ssub\dim\C Y/T(Y)$ and
$\ssub{\dim}{\RR}\ssub XT=\ssub{\dim}{\RR}\ssub YS$. By Lemma~3.1
and Theorem~2.2, $T_0$ and $S_0$ are similar. Hence there exists an
invertible linear operator $G_0:\ssub XT\to \ssub YS$ such that
\begin{equation}
\text{$G_0Tx=SG_0x$ \ for any \ $x\in \ssub XT$.} \label{g0}
\end{equation}
Pick linear subspaces $E$ and $F$ in $X$ and $Y$ respectively such
that $E\oplus T(X)=X$ and $F\oplus T(Y)=Y$. Since $\ssub\dim\C
X/T(X)=\ssub\dim\C Y/T(Y)$, we have $\ssub\dim\C E=\ssub\dim\C F$
and therefore there exists an invertible linear operator $G_1:E\to
F$. Let $\{x_\alpha\}_{\alpha\in A}$ be a maximal $T$-independent
family in $X$. According to Lemma~3.4, for any $\alpha\in A$, there
exists a unique sequence $\{x_{\alpha,n}\}_{n\in\zp}$ of elements of
$E$ such that
\begin{equation}
x_{\alpha}\equiv \sum_{j=0}^n T^j x_{\alpha,j}\ (\bmod\,
T^{n+1}(X))\ \ \text{for each}\ \ n\in\zp.\label{xalpha}
\end{equation}
Since $S$ is tame, for any $\alpha\in A$, there exists $y_\alpha\in
Y$ such that
\begin{equation}
y_{\alpha}\equiv \sum_{j=0}^n S^j G_1x_{\alpha,j}\ (\bmod\,
S^{n+1}(Y))\ \ \text{for each}\ \ n\in\zp.\label{yalpha}
\end{equation}

Let $x\in X$ and $\Lambda=\Lambda(x)$ be the finite subset of $A$
furnished by Lemma~3.3. Pick a positive integer $n$, $w_n\in \ssub
XT$ and non-zero $r_{\alpha,n}\in\RR_+$ for $\alpha\in\Lambda$ such
that (\ref{tin}) is satisfied. Let
\begin{equation}
y_n=G_0w_n+\sum_{\alpha\in\Lambda}r_{\alpha,n}(S)y_\alpha.
\label{yn}
\end{equation}
First, let us verify that $y_n\in S^n(Y)$. According to the
classical Taylor theorem, for any $\alpha\in\Lambda$, there exists a
unique polynomial
$$
p_{\alpha,n}(z)=\sum_{j=0}^{n-1} a_{\alpha,j}z^j
$$
such that $r_{\alpha,n}-p_{\alpha,n}$ has zero of order at least $n$
in zero. Then $r_{\alpha,n}(T)x_\alpha\equiv
p_{\alpha,n}(T)x_\alpha\,\, (\bmod\, T^n(X))$. Using (\ref{xalpha}),
we have
$$
r_{\alpha,n}(T)x_\alpha\equiv p_{\alpha,n}(T)x_\alpha\equiv
\sum_{j=0}^{n-1}T^j\biggl( \sum_{k=0}^j
a_{\alpha,k}x_{\alpha,j-k}\biggr)\ (\bmod\, T^n(X)).
$$
Summing over $\alpha$, we obtain
\begin{equation}
\sum_{\alpha\in\Lambda}r_{\alpha,n}(T)x_\alpha\equiv
\sum_{j=0}^{n-1}T^j\biggl( \sum_{\alpha\in\Lambda}\sum_{k=0}^j
a_{\alpha,k}x_{\alpha,j-k}\biggr)\ (\bmod\, T^n(X)). \label{for1}
\end{equation}
Similarly using (\ref{yalpha}) instead of (\ref{xalpha}), we have
\begin{equation}
\sum_{\alpha\in\Lambda}r_{\alpha,n}(S)y_\alpha\equiv
\sum_{j=0}^{n-1}S^jG_1\biggl( \sum_{\alpha\in\Lambda}\sum_{k=0}^j
a_{\alpha,k}x_{\alpha,j-k}\biggr)\ (\bmod\, S^n(Y)). \label{for2}
\end{equation}
From (\ref{tin}) and (\ref{for1}) it follows that
$$
\sum_{j=0}^{n-1}T^j\biggl( \sum_{\alpha\in\Lambda}\sum_{k=0}^j
a_{\alpha,k}x_{\alpha,j-k}\biggr)\equiv0\ (\bmod\, T^n(X)).
$$
According to the uniqueness part of Lemma~3.4,
$$
\sum_{\alpha\in\Lambda}\sum_{k=0}^j a_{\alpha,k}x_{\alpha,j-k}=0\ \
\text{for $0\leq j\leq n-1$.}
$$
Substituting this into (\ref{for2}) and using (\ref{yn}), we obtain
$y_n\equiv 0\,\,(\bmod\, S^n(Y))$, or equivalently, $y_n\in S^n(Y)$.
Thus, there exists a unique $y\in Y$ for which $S^ny=y_n$. We write
$Gx=y$.

Let us check that the map $G$ is well-defined. To this end it
suffices to verify that $y$ does not depend on $n$. Let $m>n$ and
the vectors $y_n$ and $y_m$ be defined by (\ref{yn}). According to
Lemma~3.3, $u_m=T^{m-n}u_n$ and $r_{\alpha,m}=M^{m-n}r_{\alpha,n}$
for each $\alpha\in\Lambda$. Substituting these equalities into
(\ref{yn}) and taking (\ref{g0}) into account, we see that
$y_m=S^{m-n}y_n$ and therefore $y$ does not depend on the choice of
$n$. The map $G:X\to Y$ is defined. The linearity of $G$ follows
easily from its definition. Thus, we have a linear operator $G:X\to
Y$.

First, let us show that $GT=SG$. Let $x\in X$, $\Lambda$ be the
finite set furnished by Lemma~3.3 and $n\in\zp$, $u_n\in\ssub XT$,
$r_{\alpha,n}\in\RR_+\setminus\{0\}$ for $\alpha\in\Lambda$ be such
that (\ref{tin}) is satisfied. By definition of $G$, we have $Gx=y$,
where $S^ny=y_n$ and $y_n$ is defined in (\ref{yn}). Then
$S^{n-1}SGx=y_n$. By (\ref{tin}),
$$
T^{n-1}Tx=w_n+\sum_{\alpha\in\Lambda}r_{\alpha,n}(T)x_\alpha.
$$
From the definition of  $G$ and (\ref{yn}) it follows that
$S^{n-1}GTx=y_n$. Hence $S^{n-1}GTx=S^{n-1}SGx$. Since $S$ is
injective, we have $GTx=SGx$. Thus, $GT=SG$.

Next, we shall verify the injectivity of $G$. Suppose that $x\in X$
and $Gx=0$. First, let us consider the case $x\notin \ssub XT$. In
this case the finite subset $\Lambda=\Lambda(x)$ of $A$, provided by
Lemma~3.3, is non-empty. Let $n\in\zp$, $u_n\in\ssub XT$,
$r_{\alpha,n}\in\RR_+\setminus\{0\}$ for $\alpha\in\Lambda$ be such
that (\ref{tin}) is satisfied. Since $x\notin\ssub XT$, there exists
a positive integer $m$ such that $T^nx\notin T^m(X)$. Choose
polynomials
$$
p_{\alpha,m}(z)=\sum_{j=0}^{m-1} a_{\alpha,j}z^j
$$
such that $r_{\alpha,n}-p_{\alpha,m}$ has zero of order at least $m$
in zero. As above, we have
\begin{align}
T^nx&\equiv \sum_{\alpha\in\Lambda}r_{\alpha,n}(T)x_\alpha\equiv
\sum_{j=0}^{m-1}T^j\biggl( \sum_{\alpha\in\Lambda}\sum_{k=0}^j
a_{\alpha,k}x_{\alpha,j-k}\biggr)\ (\bmod\, T^m(X))\ \ \text{and}
\label{AA}
\\
S^nGx&\equiv\sum_{\alpha\in\Lambda}r_{\alpha,n}(S)y_\alpha\equiv
\sum_{j=0}^{m-1}S^jG_1\biggl( \sum_{\alpha\in\Lambda}\sum_{k=0}^j
a_{\alpha,k}x_{\alpha,j-k}\biggr)\ (\bmod\, S^m(Y)).\label{BB}
\end{align}
Since $T^nx\notin T^m(X)$, (\ref{AA}) implies that there exists $j$,
$0\leq j\leq m-1$ such that
$$
\sum_{\alpha\in\Lambda}\sum_{k=0}^j a_{\alpha,k}x_{\alpha,j-k}\neq0.
$$
Since $G_1$ is injective,
$$
G_1\sum_{\alpha\in\Lambda}\sum_{k=0}^j
a_{\alpha,k}x_{\alpha,j-k}\neq0.
$$
Using (\ref{BB}) and the last display, we have $S^nGx\notin S^m(Y)$.
Therefore $Gx\neq 0$. This contradiction shows that $x\in \ssub XT$.
Since the restriction of $G$ to $\ssub XT$ coincides with $G_0$ and
$G_0$ is injective, we obtain $x=0$. Injectivity of $G$ is proven.

Finally let us prove that $G$ is onto. Since the restriction of $G$
to $\ssub XT$ coincides with $G_0$ and the restriction of $G$ to $E$
coincides with $G_1$, from the definition of $G$ and the equality
$GT=SG$ it follows that
\begin{equation}
\begin{array}{l}
\text{if $x\in X$ and the sequence $\{x_n\}_{n\in\zp}$ of elements
of $E$ satisfy (\ref{tame1}) for any $n\in\zp$,}\\ \text{then
$Gx\equiv G_1x_0+{\dots}+S^nG_1x_n\ (\bmod\, S^{n+1}(Y))$ for each
$n\in\zp$.}
\end{array}
\label{P}
\end{equation}
Let $y\in Y$. According to Lemma~3.4 there exists a sequence
$\{y_n\}_{n\in\zp}$ of elements of $F$ for which $y\equiv
y_0+{\dots}+S^ny_n\ (\bmod\, S^{n+1}(Y))$ for any $n\in\zp$. Since
$T$ is tame, there exists $x\in X$ such that (\ref{tame1}) is
satisfied for any $n\in\zp$ with $x_n=G_1^{-1}y_n$. Using (\ref{P}),
we see that
$$
Gx\equiv y\equiv y_0+{\dots}+S^ny_n\ (\bmod\, S^{n+1}(Y))\ \
\text{for any $n\in\zp$}.
$$
Hence $Gx-y\in \ssub YS$. Therefore $G(x-G_0^{-1}(Gx-y))=Gx-Gx+y=y$
and $y$ is in the range of $G$. Surjectivity  of $G$ is proven.
Thus, $G:X\to Y$ is an invertible linear operator and $GT=SG$ and
therefore $T$ and $S$ are similar. The proof of Theorem~3.2 is
complete.

\section{Auxiliary results}

Recall that a Fr\'echet space is a complete metrizable locally
convex topological vector space. We shall prove several lemmas,
which will be used in the proof of Theorems~1.2 and~1.5.

\subsection{Algebraic dimensions}

\lemma{4.1}Let $\{a_n\}_{n\in\zp}$ be a sequence of non-negative
numbers such that the set $\{n\in\zp:a_n>0\}$ is infinite and
$a_n^{1/n}\to 0$ as $n\to\infty$. Then there exists an infinite set
$A\subseteq\zp$ such that $a_n>0$ for each $n\in A$ and
$\sum\limits_{k=n+1}^\infty a_k=o(a_n)$ as $n\to\infty$, $n\in A$.
\epr

\proof For $\alpha>0$ denote $c_\alpha=\sup\limits_{n\in\zp}
a_n\alpha^{-n}$. Since $a_n^{1/n}\to 0$ as $n\to\infty$ and there
are positive $a_n$, we see that for any $\alpha>0$,
$c_\alpha\in(0,\infty)$ and there is $n(\alpha)\in\zp$ for which
$a_{n(\alpha)}=c_\alpha\alpha^{n(\alpha)}$.

Suppose now that the required set $A$ does not exist. Then there
exists $c>0$ such that $\sum\limits_{k=n+1}^\infty a_k\geq ca_n$ for
each $n\in\zp$. Since $a_n\leq c_\alpha\alpha^n$ for any $\alpha>0$
and any $n\in\zp$, we have
$$
ca_{n(\alpha)}=cc_\alpha \alpha^{n(\alpha)}\leq
\sum_{k=n(\alpha)+1}^\infty a_k\leq \sum_{k=n(\alpha)+1}^\infty
c_\alpha\alpha^{k+1}= \frac{\alpha
c_\alpha}{1-\alpha}\alpha^{n(\alpha)}
$$
for any $\alpha\in(0,1)$. Hence $c\leq \alpha/(1-\alpha)$ for each
$\alpha\in (0,1)$, which is impossible since $c$ is positive. This
contradiction completes the proof. \square

\lemma{4.2}Let $X$ be a Fr\'echet space and $\{x_n\}_{n\in\zp}$ be a
linearly independent sequence in $X$. Then there exists a sequence
$\{t_n\}_{n\in\zp}$ of positive numbers such that for any sequence
$\{b_n\}_{n\in\zp}$ of complex numbers,
\begin{align}
&\text{if}\ \sum_{n=0}^\infty t_n|b_n|<\infty\ \text{then the
series}\ \sum_{n=0}^\infty b_nx_n\text{\ is absolutely convergent
in}\ X;\label{A}
\\
&\text{if}\ \sum_{n=0}^\infty t_n|b_n|<\infty\ \text{and}\
\sum_{n=0}^\infty b_nx_n=0\ \text{then}\ b_n=0\text{\ for each\ \
}n\in\zp.\label{B}
\end{align}\rm

\proof Let $\{p_n\}_{n\in\zp}$ be a sequence of seminorms defining
the topology of $X$ such that $p_{n+1}(x)\geq p_n(x)$ for any $x\in
X$ and  $n\in\zp$. Since $x_n$ are linearly independent, we can
without loss of generality, assume that
\begin{align}
&\inf\{p_n(x_n-y):y\in\li\{x_1,\dots,x_{n-1}\}\}\geq\epsilon_n\in(0,1];
\label{alpha}
\\
&p_n(x_n)=1\text{\ for all}\ n\in\zp.\label{beta}
\end{align}
Indeed, if it is not the case, we can replace $p_n$ by a subsequence
$p_{k_n}$ to make (\ref{alpha}) valid and then take $x_n/p_n(x_n)$
instead of $x_n$ to make (\ref{beta}) valid.

Evidently, there exists an increasing sequence $\{t_n\}_{n\in\zp}$
of positive numbers such that
\begin{align}
&\sum_{n=0}^\infty\frac{p_k(x_n)}{t_n}<+\infty \text{\ \ for any\ \
}k\in\N; \label{gamma}
\\
&\lim_{n\to\infty}\frac{\epsilon_n t_{n+1}}{t_n}=+\infty.
\label{delta}
\end{align}

We shall show that the sequence $\{t_n\}_{n\in\zp}$ has all desired
properties. The condition (\ref{A}) follows from (\ref{gamma}). Let
us prove (\ref{B}). Suppose that $\sum\limits_{n=0}^\infty
t_n|b_n|<\infty$, $\sum\limits_{n=0}^\infty b_nx_n=0$ and there
exists $m\in\zp$ for which $b_m\neq 0$. Since $x_n$ are linearly
independent, there are infinitely many $m\in\N$ such that
$b_m\neq0$. From (\ref{delta}) and the inequality
$\sum\limits_{n=0}^\infty t_n|b_n|<\infty$ it follows that
$\lim\limits_{n\to\infty}
\Bigl(b_n\prod\limits_{k=0}^{n-1}\epsilon_k^{-1}\Bigr)^{1/n}=0$.
Applying Lemma~4.1 to the sequence $a_0=|b_0|$,
$a_n=|b_n|\epsilon_0^{-1}\cdot{\dots}\cdot\epsilon_{n-1}^{-1}$ for
$n\geq 1$, we see that there exists a strictly increasing sequence
$\{n_k\}_{k\in\zp}$ of positive integers such that $b_{n_k}\neq 0$
for each $k\in\zp$ and
$\sum\limits_{m=n_k+1}^\infty|b_m|=o(\epsilon_{n_k}b_{n_k})$ as
$k\to\infty$. Thus, there exists $r\in\zp$ for which $
\sum\limits_{n=r+1}^\infty|b_n|<\frac{\epsilon_r|b_r|}{2}$. Then
\begin{align*}
&0=p_r\biggl(\sum_{n=0}^\infty
b_nx_n\biggr)\geq|b_r|\inf\{p_r(x_r-y):y\in\li\{x_1,\dots,x_{r-1}\}\}-
\sum_{n=r+1}^\infty|b_n|p_r(x_n)\geq
\\
&\geq |b_r|\epsilon_r - \sum_{n=r+1}^\infty|b_n|p_n(x_n)\geq
|b_r|\epsilon_r -\sum_{n=r+1}^\infty|b_n|>
\frac{|b_r|\epsilon_r}{2}>0.
\end{align*}
This contradiction proves (\ref{B}). \square

Lemma~4.2 gives an alternative proof of the following well-known
fact. The classical proof goes along the same lines as the proof of
a more general Lemma~4.5 below.

\cor{4.3}Let $X$ be an infinite dimensional Fr\'echet space. Then
$\ssub\dim{\C}X\geq\ccc$. \epr

\proof Since $X$ is infinite dimensional, there exists a linearly
independent sequence $\{x_n\}_{n\in\zp}$ in $X$. By Lemma~4.2 there
exists a sequence $\{t_n\}_{n\in\zp}$ of positive numbers,
satisfying (\ref{A}) and (\ref{B}). Pick a sequence
$\{c_n\}_{n\in\zp}$ of positive numbers such that
$\sum\limits_{n=0}^\infty c_nt_n<\infty$. For any $z$ from the unit
circle $\T=\{z\in\C:|z|=1\}$ consider $y_z=\sum\limits_{n=0}^\infty
c_nz^nx_n$. The series converges absolutely in $X$ according to
(\ref{A}). Since the sequences $s_z=\{c_nz^n\}_{n\in\zp}$ are
linearly independent, (\ref{B}) implies that the family of vectors
$\{y_z\}_{z\in\T}$ is linearly independent in $X$. Since $\T$ has
cardinality $\ccc$, we have $\ssub\dim{\C}X\geq\ccc$. \square

The next proposition deals with continuous linear operators with
empty spectrum acting on Fr\'echet spaces. It worth noting that such
operators do exist. For instance, let $X$ be the space of infinitely
differentiable functions $f:[0,1]\to\C$ such that $f^{(j)}(0)=0$ for
any $j\in\zp$ endowed with the topology of uniform convergence of
all derivatives. Then $X$ is a Fr\'echet space and the Volterra
operator acts continuously on $X$ and has empty spectrum.

\prop{4.4}Let $X$ be a non-zero Fr\'echet space, $T:X\to X$ be a
continuous linear operator with empty spectrum and $\ssub X\RR$ be
$X$ considered as a linear space over $\RR$ with the multiplication
$r\cdot x=r(T)x$. Then $\ssub\dim{\C} X=\ssub\dim{\RR}\ssub X\RR$.
\epr

\proof Let $\nu=\ssub\dim{\C} X$ and $\mu=\ssub\dim{\RR}\ssub X\RR$.
Since there is no operators with empty spectrum on a finite
dimensional space, we see that $X$ is infinite dimensional. By
Corollary~4.3, $\nu\geq\ccc$. Taking into account that the algebraic
dimension of $\RR$ as a linear space over $\C$ equals $\ccc$, we see
that $\mu\cdot\ccc=\nu$. Hence $\mu\leq\nu$ and $\mu=\nu$ if
$\nu>\ccc$. It remains to verify that $\mu\geq\ccc$.

Let $x\in X\setminus \{0\}$. Since $\sigma(T)=\varnothing$, the
sequence $\{T^nx\}_{n\in\zp}$ is linearly independent. By Lemma~4.2
there exists a sequence $\{t_n\}_{n\in\zp}$ of positive numbers,
satisfying (\ref{A}) and (\ref{B}) for $x_n=T^nx$. Pick a sequence
$\{c_n\}_{n\in\zp}$ of positive numbers such that
$\sum\limits_{n=m}^\infty c_{n-m}t_n<\infty$ for each $m\in\zp$ and
$\lim\limits_{n\to\infty} \frac{c_{n+1}}{c_n}=0$.

Recall that $A\subset\T$ is called independent if
$z_1^{k_1}\cdot{\dots}\cdot z_n^{k_n}\neq 1$ for any pairwise
different $z_1,\dots,z_n\in A$ and any non-zero integers
$k_1,\dots,k_n$. Let $A\subset \T$ be an independent set of
cardinality $\ccc$ (such a set can even be chosen to be compact, see
\cite{gg}). According to the classical Kronecker theorem \cite{gg}
for any pairwise different $z_1,\dots,z_n\in A$, the sequence
$\{(z_1^k,\dots,z_n^k)\}_{k\in\zp}$ is dense in $\T^n$. Therefore
for any pairwise different $z_1,\dots,z_n\in A$ and any
$w_1,\dots,w_n\in\C$,
\begin{equation}
\slim_{k\to\infty}\biggl|\sum_{j=1}^n
w_jz_j^k\biggr|=\sum_{j=1}^n|w_j|. \label{kro}
\end{equation}
For each $z\in A$ consider the vector $y_z=\sum\limits_{n=0}^\infty
c_nz^nT^nx$. The series converges absolutely in $X$ according to
(\ref{A}) for $x_n=T^nx$. In order to prove that $\mu\geq\ccc$ it
suffices to verify that the family $\{y_z\}_{z\in A}$ is linearly
independent in $\ssub X\RR$.

Suppose the contrary. Then there exist pairwise different
$z_1,\dots,z_n\in A$ and $r_1,\dots r_n\in\RR\setminus\{0\}$ such
that $r_1(T)y_{z_1}+{\dots}+r_n(T)y_{z_n}=0$. Multiplying $r_j$ by
the least common multiple of their denominators, we see that there
are non-zero polynomials $p_1,\dots,p_n$ such that
$p_1(T)y_{z_1}+{\dots}+p_n(T)y_{z_n}=0$. Let $m=\max\limits_{1\leq
j\leq n}{\rm deg}\,p_j$. Then
\begin{equation}
p_j(z)=\sum_{l=0}^m a_{j,l}z^l\ \text{for}\ 1\leq j\leq n,\
\text{where}\ a_{j,l}\in\C\ \text{and}\ \sum_{j=1}^n |a_{j,m}|>0.
\label{poly}
\end{equation}
From the definition of $y_z$ and continuity of $T$ it follows that
$$
0=\sum_{j=1}^n p_j(T)y_{z_j}=\sum_{k=0}^\infty \alpha_k T^kx,\
\text{where}\ \alpha_k=\sum_{j=1}^n\sum_{l=0}^{\min\{m,k\}}
a_{j,l}z_j^{k-l}c_{k-l}.
$$
For $k\geq m$, we have
$$
\alpha_k=\beta_k+\gamma_k,\ \text{where}\
\beta_k=c_{k-m}\sum_{j=1}^na_{j,m}z_j^{k-m}\ \text{and}\
\gamma_k=\sum_{l=0}^m c_{k-l}\sum_{j=1}^na_{j,l}z_j^{k-l}.
$$
Clearly $\gamma_k=O(c_{k-m+1})$ as $k\to\infty$. Using (\ref{kro})
and (\ref{poly}), we have
$$
\slim\limits_{k\to\infty}\frac{|\beta_k|}{c_{k-m}}=\sum_{j=1}^n
|a_{j,m}|>0.
$$
Taking into account that $c_{k-m+1}=o(c_{k-m})$ as $k\to\infty$, we
see that there is an infinite set $\Lambda\subset\zp$ such that
$\beta_k\neq 0$ for any $k\in\Lambda$ and $\gamma_k/\beta_k\to 0$ as
$k\to\infty$, $k\in\Lambda$. It follows that
$\alpha_k=\beta_k+\gamma_k\neq 0$ for sufficiently large
$k\in\Lambda$. Since $\alpha_k=O(c_{k-m})$ as $k\to\infty$ and
$\sum\limits_{k=m}^\infty c_{k-m}t_k<\infty$, we have $\alpha_k=0$
for each $k\in\zp$ according to (\ref{B}) for $x_n=T^nx$. This
contradiction proves the linear independence of $\{y_z\}_{z\in A}$
in $\ssub X\RR$. \square

\lemma{4.5}Let $X$ be a complete metrizable topological vector
space, $Y$ be a linear subspace of $X$, which is a union of
countably many closed subsets of $X$. Suppose also that $X/Y$ is
infinite dimensional. Then $\ssub\dim\C X/Y\geq\ccc$. \epr

\proof Pick a sequence $\{B_n\}_{n\in\zp}$ of closed subsets of $X$
such that $B_n\subseteq B_{n+1}$ for any $n\in\zp$ and
$Y=\bigcup\limits_{n=0}^\infty B_n$. Since $X$ is a complete
metrizable topological vector space there exists a complete metric
$d$ on $X$ defining the topology of $X$. Indeed, the topology of any
metrizable topological vector space can be defined by a
shift-invariant metric $d$, that is $d(x,y)= d(x+u,y+u)$ for any
$x$, $y$ and $u$ and any shift-invariant metric, defining the
topology of a complete metrizable topological vector space is
complete, see \cite{shifer}. Recall that the diameter of a subset
$A$ of a metric space $(X,d)$ is
$$
\diam A=\sup_{x,y\in A} d(x,y).
$$
Let
$$
\D=\bigcup_{n=1}^\infty \D_n,\ \ \text{where}\ \
\D_n=\{\epsilon=(\epsilon_0,\dots,\epsilon_n):\epsilon_j\in\{0,1\}
\text{\ \ for\ \ }0\leq j\leq n\}.
$$
We shall construct a family $\{U_\alpha\}_{\alpha\in\D}$ of open
subsets of $X$ such that for any $n\in\zp$,
\begin{itemize}
\item[{\bf(S1)}]$\U_\alpha\subset
U_{\alpha_1,\dots,\alpha_{n-1}}$ for any $\alpha\in\D_n$ if $n\geq
1$;
\item[{\bf(S2)}]$\sum\limits_{\alpha\in\D_n}z_\alpha y_\alpha\notin B_n$
for any $y_\alpha\in U_\alpha$ and any complex numbers $z_\alpha$
such that $\sum\limits_{\alpha\in\D_n}|z_\alpha|=1$;
\item[{\bf(S3)}]$\diam U_\alpha\leq 2^{-n}$ for any $\alpha\in\D_n$.
\end{itemize}

On step 0 pick $x_0,x_1\in X$ linearly independent modulo $Y$. Since
$B_0\subseteq Y$ is closed in $X$, we can choose neighborhoods $U_0$
and $U_1$ of $x_0$ and $x_1$ respectively small enough to ensure
that, $\diam U_0\leq 1$, $\diam U_1\leq 1$ and $z_0y_0+z_1y_1\notin
B_0$ if $y_0\in\U_0$, $y_1\in\U_1$ and $|z_0|+|z_1|=1$. Thus,
conditions (S1--S3) for $n=0$ are satisfied.

Suppose now that $m$ is a positive integer and $U_\alpha$ for
$\alpha\in\bigcup\limits_{k=0}^{m-1}\D_k$, satisfying (S1--S3) for
$n<m$ are already constructed. Since $Y$ has infinite codimension in
$X$, there exist $x_\alpha\in X$ for $\alpha\in\D_m$ such that
$x_\alpha\in U_{(\alpha_1,\dots,\alpha_{m-1})}$ for each
$\alpha\in\D_m$ and the vectors $\{x_\alpha\}_{\alpha\in\D_m}$ are
linearly independent modulo $Y$ (we use the obvious fact that the
linear span of any non-empty open subset of a topological vector
space is the whole space). Since $B_m\subseteq Y$ is closed in $X$,
we can choose neighborhoods $U_\alpha$ of $x_\alpha$ for
$\alpha\in\D_m$ small enough to ensure that (S1--S3) for $n=m$ are
satisfied. The induction procedure is complete.

Conditions (S1), (S3) and completeness of $(X,d)$ imply that for any
$\alpha\in\{0,1\}^{\zp}$, $\bigcap\limits_{n=1}^\infty
U_{(\alpha_1,\dots,\alpha_n)}$ is a one-element set $\{u_\alpha\}$.
Since the cardinality of the set $\{0,1\}^{\zp}$ is $\ccc$, it
suffices to verify that the family
$\{u_\alpha\}_{\alpha\in\{0,1\}^{\zp}}$ is linearly independent
modulo $Y$. Suppose the contrary. Then there exist pairwise
different $\alpha^1,\dots,\alpha^m\in\{0,1\}^{\zp}$ and complex
numbers $z_1,\dots,z_m$ such that
$$
\sum_{j=1}^m z_ju_{\alpha^j}\in Y \ \ \text{and}\ \ \sum_{j=1}^m
|z_j|=1.
$$
Since $Y$ is the union of $B_n$, there exists $l\in\zp$ such that
$$
\sum_{j=1}^m z_ju_{\alpha^j}\in B_l.
$$
Choose $k\in\zp$, $k\geq l$ for which $\alpha^{j,k}=
(\alpha^j_0,\dots\alpha^j_k)\in \D_k$ are pairwise different for
$1\leq j\leq m$. From (S1) it follows that
$u_{\alpha^j}\in\U_{\alpha^{j,k}}$ for $1\leq j\leq m$. According to
(S2),
$$
\sum_{j=1}^m z_ju_{\alpha^j}\notin B_k\supseteq B_l.
$$
The last two displays contradict each other. Thus,
$\{u_\alpha\}_{\alpha\in\{0,1\}^{\zp}}$ is linearly independent
modulo $Y$. \square

\lemma{4.6}Let $X$ be e Fr\'echet space and $Y$ be a linear subspace
of $X$ carrying a stronger Fr\'echet space topology and
$\mu=\ssub\dim{\C}X/Y$. If $\mu$ is finite then $Y$ is closed in
$X$, if $\mu$ is infinite then $\mu\geq\ccc$. \epr

\proof First, suppose that $\mu$ is finite. Then there exists a
finite dimensional linear subspace $E$ of $X$ such that $E\oplus
Y=X$ in the algebraic sense. Consider the linear map $T$ from the
Fr\'echet space $Y\times E$ to the Fr\'echet space $X$ given by
$T(u,y)=u+y$. Since the topology of $Y$ is stronger than the one
inherited from $X$, we see that $T$ is continuous. Since $E\oplus
Y=X$, we see that $T$ is invertible. According to the Banach inverse
mapping theorem \cite{rud}, $T^{-1}$ is continuous and therefore
$Y=(T^{-1})^{-1}(Y\times\{0\})$ is closed as a pre-image of the
closed set $Y\times\{0\}$ with respect to the continuous map
$T^{-1}$.

Assume now that $\mu$ is infinite. Choose a base $\{U_n\}_{n\in\zp}$
of convex symmetric neighborhoods of zero in $Y$. Symbols $\U_n$
stand for the closures of $U_n$ in $X$. For any $n\in\zp$ let $Y_n$
be the linear span of $\U_n$. Since $Y_n=\bigcup\limits_{k=1}^\infty
k\U_n$, each $Y_n$ is the union of countably many closed subsets of
$X$. Denote also $Z=\bigcap\limits_{n=0}^\infty Y_n$. Note that if
$\{V_k\}_{k\in\zp}$ is a base of neighborhoods of zero in $X$, then
$\{V_k\cap\U_n\}_{k\in\zp}$ is a base of neighborhoods\footnote{We
do not assume neighborhoods to be open. A neighborhood of a point
$x$ of a topological space $X$ is a set containing an open set
containing $x$.} of zero of a Fr\'echet space topology on $Y_n$,
stronger than the one inherited from $X$.

{\bf Case 1:} there exists $n\in\zp$ for which the codimension of
$Y_n$ in $X$ is infinite. According to Lemma~4.5 $\ssub\dim\C
X/Y_n\geq\ccc$. Since $Y\subseteq Y_n$, we have $\mu\geq\ccc$ as
required.

{\bf Case 2:} the codimension of $Z$ in $X$ is infinite and for any
$n\in\zp$ the codimension of $Y_n$ in $X$ is finite. According to
the already proven first part of the lemma, any $Y_n$ is closed in
$X$ and therefore $Z$ is closed in $X$. Hence the infinite
dimensional Fr\'echet space $X/Z$ has algebraic dimension $\geq\ccc$
according to Corollary~4.3. Since $Y\subseteq Z$, we have
$\mu\geq\ccc$ as required.

{\bf Case 3:} the codimension of $Z$ in $X$ is finite. As in the
previous case $Z$ is closed in $X$ as the intersection of closed
linear subspaces $Y_n$. Finiteness of the codimension of $Z$ implies
that $Y_n=Z$ for sufficiently large $n$. Consider the identity
embedding $J:Y\to Z$, where $Z$ carries the Fr\'echet space topology
inherited from $X$. Since $Y_n=Z$ for sufficiently large $n$, it
follows that $J$ is almost open \cite{bonet}. Since any almost open
continuous linear map between Fr\'echet spaces is surjective and
open \cite{bonet}, we see that $J$ is onto. Thus, $Y=Z$ and
therefore $Y$ has finite codimension in $X$. This contradiction
shows that Case~3 does not occur. \square

\subsection{Density of ranges}

For a family $\{\tau_\alpha\}_{\alpha\in A}$ of topologies on a set
$X$, the symbol $\vvv\limits_{\alpha\in A}\tau_\alpha$ stands for
the topology, whose base is formed by the sets
$\bigcap\limits_{j=1}^n U_j$, where $U_j\in\tau_{\alpha_j}$ and
$\alpha_j\in A$. In other words, $\tau$ is the weakest topology
stronger than each $\tau_\alpha$.

\lemma{4.7}Let $(X_n,\tau_n)$ for $n\in\zp$ be Fr\'echet spaces such
that $X_{n+1}$ is a linear subspace of $X_n$ and
$\tau_n\bigr|_{X_{n+1}}\subseteq \tau_{n+1}$ for each $n\in\zp$. Let
also $Y=\bigcap\limits_{n=0}^\infty X_n$ be endowed with the
topology $\tau=\vvv\limits_{n=0}^\infty \tau_n\bigr|_Y$. Then
$(Y,\tau)$ is a Fr\'echet space. Moreover, if $X_{n+1}$ is
$\tau_n$-dense in $X_n$ for each $n\in\zp$, then $Y$ is
$\tau_n$-dense in $X_n$ for each $n\in\zp$. \epr

\proof The topology $\tau$ is metrizable since each $\tau_n$ is
metrizable. Let $\{x_k\}_{k\in\zp}$ be a Cauchy sequence in
$(Y,\tau)$. Since $\tau$ is stronger than the restriction to $Y$ of
any $\tau_n$, $\{x_k\}_{k\in\zp}$ is a Cauchy sequence in
$(X_n,\tau_n)$ for any $n\in\zp$. Hence for each $n\in\zp$, $x_k$
converges to $u_n\in X_n$ with respect to $\tau_n$. Since the
restriction of $\tau_n$ to $X_{n+1}$ is weaker than $\tau_{n+1}$,
$x_n$ is $\tau_n$ convergent to $u_{n+1}$ in $X_n$. The uniqueness
of a limit of a sequence in a Hausdorff topological space implies
that $u_{n+1}=u_n$ for each $n\in\zp$. Hence there exists $u\in Y$
such that $u_n=u$ for each $n\in\zp$. Since $x_k$ is
$\tau_n$-convergent to $u$ for any $n\in\zp$, we see that $x_k$ is
$\tau$-convergent to $u$. The completeness of $(Y,\tau)$ is proved.

The density part follows from the Mittag--Leffler lemma. It is also
a particular case of Lemma~3.2.2 \cite{jochen}, dealing with
projective limits of sequences of complete metric spaces. \square

\prop{4.8}Let $T$ be an injective continuous linear operator with
dense range acting on a Fr\'echet space $X$. Then $\ssub XT$ is
dense in $X$. Moreover, $\ssub XT$ carries a Fr\'echet space
topology stronger than the topology inherited from $X$, with respect
to which the restriction $T_0$ of $T$ to $\ssub XT$ is continuous.
\epr

\proof Let $\{p_k\}_{k\in\zp}$ be a sequence of seminorms defining
the initial topology $\tau_0$ on $X_0=X$. Consider the topology
$\tau_n$ on $X_n=T^n(X)$ given by the sequence of seminorms
$$
p_{n,k}(x)=\sum_{j=0}^np_k(T^{-j}(x)),\quad k\in\zp.
$$
One can easily verify that $(X_n,\tau_n)$ is a Fr\'echet space and
$\tau_n\bigr|_{X_{n+1}}\subseteq \tau_{n+1}$ for each $n\in\zp$.
Moreover, density of $T(X)$ in $X$ implies $\tau_n$-density of
$T(X_{n})=X_{n+1}$ in $X_n$ for each $n\in\zp$. Indeed, the
restriction of $T$ to $X_n$ considered as a linear operator on the
Fr\'echet space $(X_n,\tau_n)$ is similar to $T$ with the continuous
invertible operator $G=T^n:X\to X_n$, providing the similarity. By
Lemma~4.7 $\ssub XT=\bigcap\limits_{n=0}^\infty X_n$ is dense in $X$
and $\ssub XT$ with the topology $\tau=\vvv\limits_{0}^\infty\tau_n$
is a Fr\'echet space. One can easily verify that $T_0$ is
$\tau$-continuous. \square

\subsection{Weak boundedness}

\defin{3}A sequence $\{f_n\}_{n\in\zp}$ of linear functionals on a
linear space $X$ is said to be {\it weakly bounded} if there exists
a sequence $\{b_n\}_{n\in\zp}$ of positive numbers such that for any
$x\in X$ there is $c(x)>0$ for which $|f_n(x)|\leq c(x)b_n$ for each
$n\in\zp$.

\lemma{4.9}Let $\{f_n\}_{n\in\zp}$ be a sequence of continuous
linear functionals on a Fr\'echet space $X$. Then
$\{f_n\}_{n\in\zp}$ is weakly bounded if and only if there exists a
continuous seminorm $p$ on $X$, with respect to which all $f_n$ are
bounded. \epr

\proof Suppose that each $f_n$ is bounded with respect to a
continuous seminorm $p$ on $X$. Then
$$
b_n=\sup_{p(x)\leq 1}|f_n(x)|<\infty\ \ \text{for any}\ \ n\in\zp.
$$
Clearly $|f_n(x)|\leq p(x)b_n$ for any $x\in X$ and any $n\in\zp$,
which means that $\{f_n\}_{n\in\zp}$ is weakly bounded.

Suppose now that $\{f_n\}_{n\in\zp}$ is weakly bounded and
$\{b_n\}_{n\in\zp}$ is a sequence of positive integers such that
$f_n(x)=O(b_n)$ as $n\to\infty$ for each $x\in X$. Then
$$
p(x)=\sup_{n\in\zp}\frac{|f_n(x)|}{b_n}
$$
is a seminorm on $X$, with respect to which all $f_n$ are bounded.
It remains to verify that $p$ is continuous. Clearly the unit
$p$-ball $W_p=\{x\in X:p(x)\leq 1\}$ satisfies
$W_p=\bigcap\limits_{n=0}^\infty U_n$, where $U_n=\{x\in
X:|f_n(x)|\leq b_n\}$. Since $f_n$ are continuous, we see that the
sets $U_n$ are closed. Therefore $W_p$ is closed. Thus, $W_p$ is a
barrel, that is a closed convex balanced absorbing subset of $X$.
Since any Fr\'echet space is barreled \cite{bonet}, $W_p$ is a
neighborhood of zero in $X$. Hence $p$ is continuous. \square

\section{Continuous tame operators}

We start with a criterion of tameness for general linear operators.

\lemma{5.1}Let $T$ be a linear operator on a linear space $X$. Then
$T$ is tame if and only if for any sequence $\{y_n\}_{n\in\zp}$ of
elements of $X\setminus T(X)$ and any sequence $\{c_n\}_{n\in\zp}$
of complex numbers there exists $x\in X$ such that
\begin{equation}
x\equiv\sum_{j=0}^n c_jT^jy_j\ (\bmod T^{n+1}(X)) \label{tame2}
\end{equation}
for each $n\in\zp$. \epr

\proof Clearly (\ref{tame1}) with $x_j=c_jy_j$ is exactly
(\ref{tame2}) and the 'if' part of Lemma~3.4 follows. It remains to
prove the 'only if' part. Suppose that for any sequence
$\{y_n\}_{n\in\zp}$ of elements of $X\setminus T(X)$ and any
sequence $\{c_n\}_{n\in\zp}$ of complex numbers there exists $x\in
X$ such that (\ref{tame2}) is satisfied for each $n\in\zp$. We have
to prove that $T$ is tame. If $X=T(X)$ the result is trivial.
Suppose that $X\neq T(X)$ and fix $u\in X\setminus T(X)$. Let
$\{x_n\}_{n\in\zp}$ be a sequence of elements of $X$. It suffices to
construct sequences $\{y_n\}_{n\in\zp}$ of elements of $X\setminus
T(X)$ and $\{c_n\}_{n\in\zp}$ of complex numbers such that for any
$n\in\zp$ validity of (\ref{tame1}) is equivalent to validity of
(\ref{tame2}).

On step $0$ we put $y_0=x_0$, $c_0=1$ if $x_0\notin T(X)$ and
$y_0=u$, $c_0=0$ if $x_0\in T(X)$. Obviously (\ref{tame1}) is
equivalent to (\ref{tame2}) for $n=0$. Suppose now that $m$ is a
positive integer and $y_0,\dots,y_{m-1}\in X\setminus T(X)$ and
$c_0,\dots,c_{m-1}\in \C$ are such that (\ref{tame1}) is equivalent
to (\ref{tame2}) for $0\leq n\leq m-1$. To say that (\ref{tame1}) is
equivalent to (\ref{tame2}) is the same as to say that
$$
\sum_{j=0}^n T^j(x_j-c_jy_j)\in T^{n+1}(X).
$$
Since (\ref{tame1}) is equivalent to (\ref{tame2}) for $n=m-1$,
there exists $w\in X$ such that
$$
\sum_{j=0}^{m-1} T^j(x_j-c_jy_j)=T^mw.
$$
If $w+x_m\in T(X)$, we put $y_m=u$ and $c_m=0$. If $w+x_m\notin
T(X)$, we put $y_m=w+x_m$ and $c_m=1$. In any case we have
$$
\sum_{j=0}^{m} T^j(x_j-c_jy_j)=T^m(u+x_m-c_my_m)\in T^{m+1}(X).
$$
Thus, (\ref{tame1}) is equivalent to (\ref{tame2}) for $n=m$. The
inductive procedure is complete and so is the proof of the lemma.
\square

\lemma{5.2}Let $T$ be a linear operator acting on a linear space
$X$, $n\in\zp$, $Y$ be a $T$-invariant linear subspace of $X$ such
that $T^n(X)\subseteq Y$ and $S:Y\to Y$ be the restriction of $T$ to
$Y$. Then tameness of $S$ implies tameness of $T$. \epr

\proof Suppose that $S$ is tame. For $n=0$ we have $S=T$ and the
result is trivial. Thus, we can assume that $n>0$. Let
$\{x_k\}_{k\in\zp}$ be a sequence of elements of $X$. Then
$\{T^nx_{n+k}\}_{k\in\zp}$ is a sequence of elements of
$T^n(X)\subseteq Y$. Since $S$ is tame, there exists $y\in Y$ such
that
$$
y-\sum_{j=0}^mT^{n+j}x_{n+j}\in S^{m+1}(Y)\subseteq T^{m+1}(X) \ \
\text{for each $m\in\zp$.}
$$
Hence,
$$
y-\sum_{l=n}^m T^{l}x_l\in T^{m+1}(X)\ \ \text{for $m\geq n$.}
$$
Let $x=y+\sum\limits_{j=0}^{n-1}T^jx_j$. From the last display it
follows that
$$
x-\sum_{j=0}^m T^jx_j\in T^{m+1}(X)\ \ \text{for any $m\in\zp$.}
$$
Hence $T$ is tame. \square

\lemma{5.3}Let $T$ be an injective linear operator on a linear space
$X$, $x\in X$ and $y_0,\dots,y_n\in X\setminus T(X)$ be such that
$(\ref{tame2})$ is satisfied for some $c_0,\dots,c_n\in\C$. Then the
numbers $c_0,\dots,c_n$ are uniquely determined by
$x,y_0,\dots,y_n$. \epr

\proof Suppose that (\ref{tame2}) is also satisfied with $c_j$
replaced by $c'_j\in\C$. We have to prove that $c_j=c'_j$ for $0\leq
j\leq n$. Since $\sum\limits_{j=0}^n (c_j-c'_j)T^jy_j\in
T^{n+1}(X)$,  we see that $(c_0-c'_0)y_0\in T(X)$. Then $c_0=c'_0$
because $y_0\notin T(X)$. Hence $\sum\limits_{j=1}^n
(c_j-c'_j)T^jy_j\in T^{n+1}(X)$ and therefore $(c_1-c'_1)Ty_1\in
T^2(X)$. Since $T$ is injective, $(c_1-c'_1)y_1\in T(X)$. Then
$c_1=c'_1$ since $y_1\notin T(X)$. Proceeding in the same way, we
obtain that $c_j=c'_j$ for each $j\leq n$. \square

We shall introduce some additional notation. Let $T$ be an injective
linear operator on a linear space $X$ and $\kappa=\{y_n\}_{n\in\zp}$
be a sequence of elements of $X\setminus T(X)$. Symbol $\ssub
XT(\kappa)$ stands for the set of $x\in X$ such that for any
positive integer $n\in\zp$ there exist $c_0,\dots,c_n\in\C$ for
which (\ref{tame2}) is satisfied. According to Lemma~5.3 the numbers
$c_j$ depend only on $x$. Thus, the maps $x\mapsto c_j$ are
well-defined linear functionals on the linear space $\ssub
XT(\kappa)$. We denote them by symbols $\Phi_j=\Phi^{T,\kappa}_j$.
In this new notation (\ref{tame2}) can be rewritten as
\begin{equation}
x\equiv\sum_{j=0}^n \Phi_j(x)T^jy_j\ (\bmod\,T^{n+1}(X)),
\label{tame3}
\end{equation}
which holds true for each $x\in \ssub XT(\kappa)$ and each
$n\in\zp$. Then Lemma~5.1 immediately implies the following
corollary.

\cor{5.4}Let $T$ be an injective linear operator on a linear space
$X$. Then $T$ is tame if and only if for any sequence
$\kappa=\{y_n\}_{n\in\zp}$ of elements of $X\setminus T(X)$ and any
sequence $\{c_n\}_{n\in\zp}$ of complex numbers, there exists
$x\in\ssub XT(\kappa)$ such that $\Phi_n^{T,\kappa}(x)=c_n$ for each
$n\in\zp$. \epr

\subsection{Tameness of continuous linear operators on Fr\'echet spaces}

In order to show that certain continuous linear operators on
Fr\'echet spaces are tame we need  an old result of Eidelheit
\cite{eid}, related to the abstract moment problem. We present it in
a slightly different form obviously equivalent to the original one.
A different proof can be found in \cite{shk}.

\prop{5.5} Let $X$ be a Fr\'echet space and $\{\phi_n\}_{n\in\zp}$
be a sequence of continuous linear functionals on $X$. Then the
following conditions are equivalent:
\begin{itemize}
\item[{\bf (E1)}]for any sequence $\{c_n\}_{n\in\zp}$ of complex numbers
there exists $x\in X$ such that $\phi_n(x)=c_n$ for each $n\in\zp;$
\item[{\bf (E2)}]the family $\{\phi_n\}_{n\in\zp}$ of functionals is
linearly independent and for any continuous seminorm $p$ on $X$ the
space
$$
E_p=\{\phi\in\li\{\phi_j:j\in\zp\}:\phi\ \text{is $p$-bounded}\}
$$
is finite dimensional.
\end{itemize}\rm

\lemma{5.6}Let $X$ be a Fr\'echet space and $T:X\to X$ be an
injective continuous linear operator such that $\ssub XT$ is dense
in $X$. Then $T$ is tame. \epr

\proof  Let $\kappa=\{y_n\}_{n\in\zp}$ be a sequence of elements of
$X\setminus T(X)$. Since $y_j\notin T(X)$, the vectors
$y_0,Ty_1,\dots,T^ny_n$ are linearly independent modulo $T^{n+1}(X)$
for any $n\in\zp$. Therefore the space
$$
E^\kappa_n=\li\{T^jy_j:0\leq j\leq n\}
$$
is $n+1$-dimensional and $E^\kappa_n\cap T^{n+1}(X)=\{0\}$. Denote
$X^\kappa_n=E^\kappa_n\oplus T^{n+1}(X)$. Clearly
$$
\ssub XT(\kappa)=\bigcap_{n=0}^\infty X^\kappa_n.
$$
Let $\{p_k\}_{k\in\zp}$ be a sequence of seminorms defining the
topology of $X$. Since $T$ is injective, any $x\in X^\kappa_n$ can
be uniquely written in the form
$$
x=\sum_{j=0}^n c_jT^jx_j+T^{n+1}u,\ \ \text{where}\ \ c_j\in\C,\
u\in X.
$$
This allows us to define seminorms
$$
p_{n,k}(x)=\sum_{j=0}^n |c_j|+p_k(u)
$$
on $X^\kappa_n$. The sequence $\{p_{n,k}\}_{k\in\zp}$ of seminorms
on $X^\kappa_n$ defines a metrizable locally convex topology
$\tau_n$ on $X^\kappa_n$. Since $T$ is continuous, it follows that
the restriction of $\tau_n$ to $X^\kappa_{n+1}$ is weaker than
$\tau_{n+1}$. Moreover, the map $x\mapsto (c_0,\dots,c_n,u)$ is an
isomorphism between $(X^\kappa_n,\tau_n)$ and the Fr\'echet space
$\C^{n+1}\times X$. Hence $(X^\kappa_n,\tau_n)$ is a Fr\'echet space
for each $n\in\zp$. By Lemma~4.7 $\ssub XT(\kappa)$ endowed with the
topology $\tau=\vvv\limits_{n=0}^\infty \tau_n$ is a Fr\'echet
space. Clearly $\tau$ is defined by the family of seminorms
$\{p_{n,k}\}_{n,k\in\zp}$. Since for any $n\in\zp$, the functional
$\Phi_n=\Phi_n^{T,\kappa}$ is bounded with respect to $p_{n,0}$, we
see that all $\Phi_n$ are $\tau$-continuous linear functionals on
$\ssub XT(\kappa)$. According to Corollary~5.4 it suffices to verify
that for any sequence $\{c_n\}_{n\in\zp}$ of complex numbers, there
exists $x\in\ssub XT(\kappa)$ for which $\Phi_n(x)=c_n$ for each
$n\in\zp$. Since $\{(T^ny_n,\Phi_n)\}_{n\in\zp}$ is a biorthogonal
sequence, $\Phi_n$ are linearly independent. Let
$E=\li\{\Phi_n:n\in\zp\}$. By Proposition~5.5 it suffices to show
that for any $\tau$-continuous seminorm $p$ on $\ssub XT(\kappa)$,
$$
\text{the space}\ \ E_p=\{\phi\in E:\text{$\phi$ is $p$-bounded}\}\
\ \text{is finite dimensional.}
$$
Let $p$ be a $\tau$-continuous seminorm on $\ssub XT(\kappa)$. Since
$\tau=\vvv\limits_{k=0}^\infty \tau_k$ and restriction of $\tau_k$
to $X^\kappa_{k+1}$ is weaker than $\tau_{k+1}$ for each $k\in\zp$,
we see that there exists $n\in\zp$ for which $p$ is
$\tau_n$-continuous.

Since $T^{n+1}$ acting from $X$ to the subspace $T^{n+1}(X)$ of
$X^\kappa_n$, endowed with the topology $\tau_n$, is an isomorphism
of Fr\'echet spaces, mapping $\ssub XT$ onto itself, we have that
$\ssub XT$ is $\tau_n$-dense in $T^{n+1}(X)$. Since any $\phi\in E$
vanishes on $\ssub XT$, we see that any $\tau_n$-continuous $\phi\in
E$ vanishes on $T^{n+1}(X)\cap \ssub XT(\kappa)$. Therefore
$$
E_p\subseteq \{\phi\in E:\text{$\phi$ is
$\tau_n$-continuous}\}\subseteq \{\phi\in E:\text{\ $\phi$ vanishes
on $T^{n+1}(X)\cap \ssub XT(\kappa)\}$}.
$$
The dimension of the last space does not exceed the codimension of
$T^{n+1}(X)\cap\ssub XT(\kappa)$ in $\ssub XT(\kappa)$, which does
not exceed the codimension of $T^{n+1}(X)$ in $X^\kappa_n$, which is
finite. \square

\prop{5.7}Let $X$ be a Fr\'echet space and $T:X\to X$ be an
injective continuous linear operator such that there exists
$m\in\zp$ for which $\overline{T^m(X)}=\overline{T^{m+1}(X)}$. Then
$T$ is tame. \epr

\proof Clearly $Y=\overline{T^m(X)}$ is a closed $T$-invariant
subspace of $X$. Let $S:Y\to Y$ be the restriction of $T$ to $Y$.
Condition $\overline{T^m(X)}=\overline{T^{m+1}(X)}$ implies that the
range of $S$ is dense. By Proposition~4.8 $\ssub YS$ is dense in
$Y$. Applying Lemma~5.6, we see that $S$ is tame. From Lemma~5.2 it
follows that $T$ is tame. \square

\subsection{Bounded tame operators on Banach spaces}

In the Banach space setting the sufficient condition of tameness in
Proposition~5.7 turns out to be also necessary.

\theorem{5.8}Let $T$ be an injective bounded linear operator on a
Banach space $X$. Then $T$ is tame if and only if there exists
$m\in\zp$ for which $\overline{T^m(X)}=\overline{T^{m+1}(X)}$. \epr

\proof By Proposition~5.7 the existence of $m\in\zp$ for which
$\overline{T^m(X)}=\overline{T^{m+1}(X)}$ implies tameness of $T$.
Suppose now that $\overline{T^m(X)}\neq\overline{T^{m+1}(X)}$ for
any $m\in\zp$. Then we can choose a sequence
$\kappa=\{y_n\}_{n\in\zp}$ in $X$ such that $T^ny_n\notin
\overline{T^{n+1}(X)}$ for each $n\in\zp$. In particular, $\kappa$
is a sequence of elements of $X\setminus T(X)$. Using the
Hahn--Banach theorem, for any $n\in\zp$,  we can find a continuous
linear functional $\phi_n$ on $X$ such that $\phi_n(T^ny_n)=1$,
$\phi_n(T^jy_j)=0$ for $j<n$ and $\phi_n(x)=0$ for each $x\in
\overline{T^{n+1}(X)}$. From the definition of the functionals
$\Phi_n=\Phi_n^{T,\kappa}$ on the Fr\'echet space $\ssub XT(\kappa)$
it follows that each $\Phi_n$ is the restriction of $\phi_n$ to
$\ssub XT(\kappa)$. Hence,
$$
|\Phi_n(x)|=|\phi_n(x)|\leq \|\phi_n\|\cdot\ssub{\|x\|}{X}\ \
\text{for any $x\in \ssub XT(\kappa)$ and any $n\in\zp$}.
$$
Therefore there exists no $x\in\ssub XT(\kappa)$ such that
$\Phi_n(x)=n\|\phi_n\|$ for each $n\in\zp$. According to
Corollary~5.4, $T$ is not tame. \square

\section{Similarity to the Volterra operator}

Let $\V$ be the class of injective quasinilpotent continuous linear
operators $T$ acting on a Fr\'echet space $X$ of algebraic dimension
$\ccc$ such that the codimension of $T(X)$ in $X$ is infinite and
there exists $m\in\zp$ for which
$\overline{T^{m+1}(X)}=\overline{T^{m}(X)}$.

\prop{6.1}Any $T\in\V$ is similar to the Volterra operator $V$
acting on the Banach space ${\cal C}=C[0,1]$. \epr

\proof Let $T\in\V$. Since the codimension of $T(X)$ in $X$ is
infinite and $T(X)$ carries a stronger Fr\'echet topology (namely,
the one transferred from $X$ by the operator $T$), from Lemma~4.6 it
follows that $\ssub\dim\C X/T(X)\geq\ccc$. Since $\ssub\dim\C
X=\ccc$, we have $\ssub\dim\C X/T(X)=\ccc$. Condition
$\overline{T^{m+1}(X)}=\overline{T^{m}(X)}$ together with
Proposition~4.8 imply that $\ssub XT$ is dense in
$\overline{T^{m}(X)}$ and therefore $\ssub XT$ is non-zero.
Moreover, $\ssub XT$ carries a stronger Fr\'echet space topology,
with respect to which the restriction $T_0$ of $T$ to $\ssub XT$ is
continuous. According to Lemma~3.1 $\sigma(T_0)=\varnothing$. By
Proposition~4.4 and Corollary~4.3, $\ssub\dim\RR\ssub XT\geq\ccc$.
Since $\ssub\dim\C X=\ccc$, we have $\ssub\dim\RR\ssub XT=\ccc$.
Using Proposition~5.7, we see that $T$ is tame.

Now from Theorem~3.2 it follows that any two operators from $\V$ are
similar. Since $V$ has non-closed range, Lemma~4.6 implies that
$V({\cal C})$ has infinite codimension in $\cal C$. Clearly $V$ is
injective, quasinilpotent and $\overline{V({\cal
C})}=\overline{V^2({\cal C})}$. Hence $V\in\V$. \square

\subsection{Proof of Theorem~1.2}

Since the range of any injective quasinilpotent operator on a Banach
space is non-closed, Lemma~4.6 implies that
$\ssub{\dim}{\C}X/T(X)\geq\ccc$. Using Proposition~6.1, we obtain
that (C2) implies (C1). Suppose that (C1) is satisfied. Since by
Theorem~5.8 $V$ acting on $C[0,1]$ is tame and tameness is a
similarity invariant, we observe that $T$ is tame. Since
quasinilpotence and injectivity are also preserved by similarity, we
see that $T$ is injective and quasinilpotent. Since $T$ is tame we,
using Theorem~5.8 once again, obtain that there exists $m\in\zp$ for
which $\overline{T^{m+1}(X)}=\overline{T^{m}(X)}$. Thus, (C2) is
satisfied. The proof of Theorem~1.2 is complete.

\section{Proof of Theorem~1.5}

We start with the following general fact.

\lemma{7.1}Let $\{H_n\}_{n\in\zp}$ be a strictly decreasing sequence
of closed finite codimensional linear subspaces of a Fr\'echet space
$X$ such that $\bigcap\limits_{n=0}^\infty H_n=\{0\}$. Let also $Y$
be a Fr\'echet space and $G:X\to Y$ be an invertible linear operator
such that $G(H_n)$ is closed in $Y$ for any $n\in\zp$. Then $G$ is
continuous. \epr

\proof Choose a decreasing sequence $\{K_m\}_{m\in\zp}$ of closed
linear subspaces of $X$ such that $K_0=X$, $\ssub{\dim}{\C}
K_m/K_{m+1}=1$ for each $n\in\zp$ and $\{H_n\}_{n\in\zp}$ is a
subsequence of $\{K_m\}_{m\in\zp}$: $H_n=K_{m_n}$ for some strictly
increasing sequence of non-negative integers $\{m_n\}_{n\in\zp}$.
Pick $x_n\in K_n\setminus K_{n+1}$. Using the Hahn--Banach theorem,
we can choose continuous linear functionals $\phi_n:X\to\C$ such
that $\phi_n(x_n)=1$, $\phi_n(x_j)=0$ for $j<n$ and
$\phi_n\bigr|_{K_{n+1}}\equiv0$. Consider the functionals
$\psi_n:Y\to\C$ defined by the formulas $\psi_n(y)=\phi_n(G^{-1}y)$.
Each $\psi_n$ vanishes on a finite codimensional closed linear space
of the shape $G(H_{m})$ and therefore is continuous. Since
$$
\bigcap_{m=0}^\infty {\rm ker}\,\phi_m=\bigcap_{m=0}^\infty K_m=
\bigcap_{n=0}^\infty H_n=\{0\},
$$
we see that the functionals $\{\phi_n\}_{n\in\zp}$ separate points
of $X$. Since $G$ is invertible, it follows that
$\{\psi_n\}_{n\in\zp}$ separate points of $Y$. Consider the
topologies $\sigma$ on $X$ and $\sigma^*$ on $Y$ defined by the
families of functionals $\{\phi_n\}_{n\in\zp}$ and
$\{\psi_n\}_{n\in\zp}$ respectively. Recall that this means that
$\sigma$ and $\sigma^*$ are given by the seminorms
$$
p_k(x)=\sum_{j=0}^k |\phi_j(x)|\ \ \text{and}\ \
q_k(y)=\sum_{j=0}^k|\psi_j(y)|\ \ \ (k\in\zp)
$$
respectively. Continuity of $\phi_n$ and $\psi_n$ implies that
$\sigma$ and $\sigma^*$ are weaker than the initial topologies on
$X$ and $Y$. Since $\phi_n$ separate the points of $X$ and $\psi_n$
separate the points of $Y$, the topologies $\sigma$ and $\sigma^*$
are Hausdorff. From the formula $\psi_n(y)=\phi_n(G^{-1}y)$ it
immediately follows that $G$ is $\sigma$-$\sigma^*$ continuous.
Therefore the graph $\Gamma_G$ is closed in the product of
$(X,\sigma)$ and $(Y,\sigma^*)$ as is the graph of any continuous
map between Hausdorff topological spaces. Since $\sigma$ and
$\sigma^*$ are weaker than the initial topologies on $X$ and $Y$, we
see that $\Gamma_G$ is closed in $X\times Y$. The Banach Closed
Graph Theorem \cite{rud,bonet} implies that $G$ is continuous.
\square

We are going to use notation from Section~5. For an injective linear
operator $T$ on a linear space $X$ symbol $B(T,X)$ stands for the
set of vectors $y\in X$ such that either $y\in T(X)$ or $y\notin
T(X)$ and there exists a sequence $\{u_n\}_{n\in\zp}$ of elements of
$X$ for which the sequence $\{\Phi_n^{T,\kappa}\}_{n\in\zp}$ of
linear functionals on the space $\ssub XT(\kappa)$ is not weakly
bounded, where $\kappa=\{y+Tu_n\}_{n\in\zp}$.

For a bounded linear operator $T$ on a Banach space $X$ we denote
$$
A(T,X)= \bigcup_{n=0}^\infty T^{-n}(\overline{T^{n+1}(X)}).
$$
Clearly $A(T,X)$ is a linear subspace of $X$ as a union of the
increasing sequence of linear subspaces
$T^{-n}(\overline{T^{n+1}(X)})$.

\lemma{7.2}Let $T$ be an injective bounded linear operator on a
Banach space $X$. Then $B(T,X)\subseteq A(T,X)$. \epr

\proof Let $y\in X\setminus A(T,X)$. Then
$T^ny\notin\overline{T^{n+1}(X)}$ for each $n\in\zp$. Let also
$\{u_n\}_{n\in\zp}$ be a sequence of elements of $X$ and
$y_n=y+Tu_n$. Then $T^ny_n\notin \overline{T^{n+1}(X)}$ for each
$n\in\zp$. Using the Hahn--Banach theorem, we can choose continuous
linear functionals $\phi_n$ on $X$ such that $\phi_n(T^ny_n)=1$,
$\phi_n(T^jy_j)=0$ for $j<n$ and
$\phi_n\bigr|_{\overline{T^{n+1}(X)}}\equiv 0$. Then the functionals
$\Phi_n=\Phi_n^{T,\kappa}$ on $\ssub XT(\kappa)$ for
$\kappa=\{y_n\}_{n\in\zp}$ coincide with the restrictions of
$\phi_n$ to $\ssub XT(\kappa)$. Hence $|\Phi_n(x)|=|\phi_n(x)|\leq
\|\phi_n\|\|x\|$ for any $x\in\ssub XT(\kappa)$ and any $n\in\zp$
and therefore the sequence $\{\Phi_n\}_{n\in\zp}$ is weakly bounded.
Thus, $y\notin B(T,X)$, which proves the desired inclusion. \square

\lemma{7.3}Let $T$ be an injective bounded linear operator on a
Banach space $X$ and $\kappa=\{y_n\}_{n\in\zp}$ be a sequence of
elements of $X\setminus T(X)$ such that the sequence
$\{\Phi_n=\Phi_n^{T,\kappa}\}_{n\in\zp}$ of linear functionals on
$\ssub XT(\kappa)$ is weakly bounded. Then there exists $k\in\zp$
such that for any $n\geq k$, $T^{n-k}y_n$ does not belong to the
closure of \ $\li\{T^{n-k+j}y_{n+j}:j\geq 1\}$ \ in $X$. \epr

\proof As in Section~5, $E_n^\kappa=\li\{T^jy_j:0\leq j\leq n\}$ and
$X_n^\kappa=E_n^\kappa\oplus T^{n+1}(X)$. We endow $X_n^\kappa$ with
the norm
$$
p_n\biggl(T^{n+1}u+\sum_{j=0}^n c_jT^jy_j\biggr)=\ssub{\|u\|}X+
\sum_{j=0}^n |c_j|.
$$
Then $(X_n^\kappa,p_n)$ is a Banach space. Indeed, the map
$(c,u)\mapsto T^{n+1}u+\sum\limits_{j=0}^n c_jT^jy_j$ from
$\C^{n+1}\times X$ is an isomorphism of normed spaces and
$\C^{n+1}\times X$ is complete. Continuity of $T$ implies that there
exist $c_n>0$ for which
\begin{equation}
p_n(x)\leq c_n p_{n+1}(x)\ \ \text{for any $n\in\zp$ and any $x\in
X_{n+1}^\kappa$}.\label{pn}
\end{equation}
According to Lemma~4.7 $\ssub XT(\kappa)=\bigcap\limits_{n=0}^\infty
X_n^\kappa$ endowed with the topology defined by the sequence
$\{p_n\}_{n\in\zp}$ of norms is a Fr\'echet space. The functionals
$\Phi_n=\Phi_n^{T,\kappa}:\ssub XT(\kappa)\to\C$ are
$\tau$-continuous since each $\Phi_n$ is $p_n$-bounded. Since the
sequence $\{\Phi_n\}_{n\in\zp}$ is weakly bounded, Lemma~4.9 implies
the existence of a $\tau$-continuous seminorm $p$ on $\ssub
XT(\kappa)$, with respect to which each $\Phi_n$ is bounded.
According to (\ref{pn}), there is a positive integer $k$ such that
each $\Phi_n$ is $p_{k-1}$-bounded. Hence for any $n\in\zp$, ${\rm
ker}\,\Phi_n$ is $p_{k-1}$-closed in $\ssub XT(\kappa)$. Since
$$
\text{$T^ny_n\notin {\rm ker}\,\Phi_n$ \ and \ $\li\{T^my_m:m\neq
n+1\}\subset {\rm ker}\,\Phi_n$},
$$
we see that $T^ny_n$ does not belong to the $p_{k-1}$-closure of
$\li\{T^my_m:m\geq n+1\}$. From the definition of the norm $p_{k-1}$
it follows that for $n\geq k$, the last condition means exactly that
$T^{n-k}y_n$ does not belong to the closure in $X$ of
$$
\li\{T^{m-k}y_{m}:m\geq n+1\}=\li\{T^{n-k+j}y_{n+j}:j\geq 1\}.\ \
\square
$$

Let $T$ be an injective bounded linear operator on a Banach space
$X$. It is easy to see that for any $n\in\zp$, the operator
$$
T_n:X/T^{-n}(\overline{T^{n+1}(X)})\to \overline{T^{n}(X)}/
\overline{T^{n+1}(X)}, \quad T_n(x+T^{-n}(\overline{T^{n+1}(X)}))=
T^nx+\overline{T^{n+1}(X)}
$$
is an injective bounded linear operator with dense range. It follows
that if one of the spaces $X/T^{-n}(\overline{T^{n+1}(X)})$ or
$\overline{T^{n}(X)}/ \overline{T^{n+1}(X)}$ is finite dimensional
then the other has the same dimension and $T_n$ is an isomorphism
between them. Moreover for each $n\in\zp$, the operator
$$
\widetilde T_n:\overline{T^{n}(X)}/ \overline{T^{n+1}(X)}\to
\overline{T^{n+1}(X)}/ \overline{T^{n+2}(X)}, \quad \widetilde
T_n(x+\overline{T^{n+1}(X)})=Tx+\overline{T^{n+2}(X)}
$$
is a bounded linear operator with dense range. Thus, if
$\overline{T^{n}(X)}/ \overline{T^{n+1}(X)}$ is finite dimensional
then so is $\overline{T^{n+1}(X)}/ \overline{T^{n+2}(X)}$ and
$\widetilde T_n$ is onto. In particular, the dimension of
$\overline{T^{n+1}(X)}/ \overline{T^{n+2}(X)}$ does not exceed the
dimension of $\overline{T^{n}(X)}/ \overline{T^{n+1}(X)}$ and these
dimensions coincide if and only if $\widetilde T_n$ is an
isomorphism. In this case
$T^{-n}(\overline{T^{n+1}(X)})=T^{-n-1}(\overline{T^{n+2}(X)})$.
These observations are summarized in the following lemma.

\lemma{7.4}Let $m\in\zp$ and $T$ be an injective bounded linear
operator on a Banach space $X$. For $n\in\zp$ let
$\delta_n=\ssub\dim{\C}\overline{T^{n+1}(X)}/\overline{T^{n}(X)}$
and $\delta'_n=\ssub\dim{\C}X/T^{-n}(\overline{T^{n+1}(X)})$.
Suppose also $\min\{\delta_m,\delta'_m\}$ is finite. Then
$\delta_n=\delta'_n$ and $\delta_{n+1}\leq\delta_n$ for any $n\geq
m$. In particular, there exists $m_0\in\zp$ and $k_0\in\zp$ such
that $\delta_n=k_0$ for each $n\geq m_0$. Moreover, $k_0$ is exactly
the codimension of $A(T,X)$ in $X$ and if $n\geq m_0$ and $y\in
\overline{T^{n+1}(X)}\setminus\overline{T^{n}(X)}$, then $Ty\in
\overline{T^{n+2}(X)}\setminus\overline{T^{n+1}(X)}$. \epr

\lemma{7.5}Let $T$ be an injective bounded linear operator on a
Banach space $X$ such that $(\ref{151})$ and $(\ref{152})$ are
satisfied. Consider $H_n\subseteq X$ defined inductively by the
formulas $H_0=X$, $H_{k+1}=A(T,H_k)=\bigcup\limits_{n=0}^\infty
T^{-n}(\overline{T^{n+1}(H_k)})$ for $k\in\zp$. Then
$\{H_n\}_{n\in\zp}$ is a decreasing sequence of finite codimensional
closed linear subspaces of $X$ and $\bigcap\limits_{k=0}^\infty
H_k=\{0\}$. \epr

\proof Clearly $H_1$ is a linear subspace of $X$ as a union of an
increasing sequence $T^{-n}(\overline{T^{n+1}(X)})$ of linear
subspaces. From (\ref{152}) and Lemma~7.4 the sequence
$T^{-n}(\overline{T^{n+1}(X)})$ stabilizes and
$T^{-n}(\overline{T^{n+1}(X)})$ has finite codimension for
sufficiently large $n$. Thus, $A(T,X)=T^{-n}(\overline{T^{n+1}(X)})$
for sufficiently large $n$ and $H_1=A(T,X)$ is a closed finite
codimensional subspace of $X$. Applying the same argument
consecutively to the restrictions\footnote{We use the obvious
observation that if $T$ satisfies (\ref{151}) and (\ref{152}) then
so does any restriction of $T$ to any invariant closed finite
codimensional subspace.} of $T$ to $H_k$, $k=1,2,\dots$ we see that
$H_{k+1}$ is closed and finite codimensional in $H_k$ for each
$k\in\zp$. Thus, $\{H_n\}_{n\in\zp}$ is a decreasing sequence of
finite codimensional closed linear subspaces of $X$. It remains to
show that $\bigcap\limits_{k=0}^\infty H_k=\{0\}$.

Using Lemma~7.4, we see that for any fixed $k\in\zp$ there exists
$m_k\in\zp$ such that $T^{-n}(\overline{T^{n+1}(H_k)})$ does not
depend on $n$ provided $n\geq m_k$. Thus,
$H_{k+1}=T^{-n}(\overline{T^{n+1}(H_k)})$ for $n\geq m_k$. From this
equality and the definition of $H_k$ it follows that for any
$k\in\zp$,
$$
H_k\subseteq T^{-n}(\overline{T^{n+k}(X)})\ \ \text{for sufficiently
large $n$.}
$$
Let $y\in X$, $y\neq 0$. Using (\ref{151}) we see that there exists
$m\geq m_0$ such that $T^{m_0}y\in
\overline{T^{m}(X)}\setminus\overline{T^{m+1}(X)}$. From the last
statement in Lemma~7.4 it follows that $T^{m_0+l}y\in
\overline{T^{m+l}(X)}\setminus\overline{T^{m+l+1}(X)}$ for each
$l\in\zp$. In particular, $y\notin
T^{-m_0-l}(\overline{T^{m+l+1}(X)})$ for each $l\in\zp$. According
to the last display this means that $y\notin H_k$ for $k>m-m_0$.
Thus, $\bigcap\limits_{k=0}^\infty H_k=\{0\}$. \square

\lemma{7.6}Let $m\in\zp$, $T$ be an injective bounded linear
operator on a Banach space $X$ such that the codimension of
$\overline{T^{m+1}(X)}$ in $\overline{T^m(X)}$ is finite, $y\in
T^{-m}(\overline{T^{m+1}(X)})$ and $A$ be an infinite subset of
$\zp$ such that $0\in A$. Then there exists a sequence
$\{u_n\}_{n\in A}$ of elements of $X$ such that $T^m(y+Tu_0)$
belongs to the closed linear span of $\{T^{m+n}(y+Tu_n):n\in
A\setminus\{0\}\}$. \epr

\proof First, let us show that
\begin{equation}
T^n(X)+\overline{T^k(X)}=\overline{T^n(X)}\ \ \text{for any}\ \
k,n\in\zp,\ k\geq n\geq m. \label{tnm}
\end{equation}
The inclusion $T^n(X)+\overline{T^k(X)}\subseteq \overline{T^n(X)}$
is obvious. According to Lemma~7.4, $\overline{T^k(X)}$ has finite
codimension in $\overline{T^n(X)}$. Therefore
$T^n(X)+\overline{T^k(X)}$ is a linear subspace of
$\overline{T^n(X)}$, containing a closed finite codimensional
subspace. Hence $T^n(X)+\overline{T^k(X)}$ is closed in
$\overline{T^n(X)}$. On the other hand $T^n(X)+\overline{T^k(X)}$
contains $T^n(X)$ and therefore is dense in $\overline{T^n(X)}$.
Thus, (\ref{tnm}) is satisfied.

Let $\{j_n\}_{n\in\zp}$ be a strictly increasing sequence of
non-negative integers such that $A=\{j_n:n\in\zp\}$. Clearly
$j_0=0$.

We shall construct inductively a sequence $\{u_{j_n}\}_{n\in\zp}$ of
elements of $X$ such that for any $n\in\zp$,
\begin{equation}
\sum\limits_{k=0}^n T^{j_k+m}(y+Tu_{j_k})\in
\overline{T^{m+j_{n+1}+1}(X)}\ \ \text{and}\ \ \biggl\|
\sum\limits_{k=0}^n T^{j_k+m}(y+Tu_{j_k}) \biggr\|\leq 2^{-n}.
\label{indu}
\end{equation}

Since $y\in T^{-m}(\overline{T^{m+1}(X)})$, we have $T^my\in
\overline{T^{m+1}(X)}$. According to (\ref{tnm}), $T^ny\in
T^{m+1}(X)+\overline{T^{m+j_1+1}(X)}$. Therefore there exists
$w_0\in X$ for which $T^my+T^{m+1}w_0\in \overline{T^{m+j_1+1}(X)}$.
Now we can choose $v_0\in X$ such that
$\|T^my+T^{m+1}w_0-T^{m+j_1+1}v_0\|\leq 1$. Denoting
$u_{j_0}=w_0-T^{j_1}v_0$, we obtain $T^m(y+Tu_{j_0})\in
\overline{T^{m+j_1+1}(X)}$ and $\|T^m(y+Tu_{j_0})\|\leq 1$, that is
(\ref{indu}) for $n=0$ is satisfied. The basis of induction is
constructed. Suppose now that $q$ is a positive integer and
$u_{j_0},\dots,u_{j_{q-1}}$, satisfying (\ref{indu}) for $n\leq q-1$
are already constructed. Denote $x=\sum\limits_{k=0}^{q-1}
T^{j_k+m}(y+Tu_{j_k})$. From (\ref{indu}) for $n=q-1$ it follows
that $x\in \overline{T^{m+j_q+1}(X)}$. According to (\ref{tnm}),
$x\in T^{m+j_q+1}(X)+\overline{T^{m+j_{q+1}+1}(X)}$. Since $T^my\in
\overline{T^{m+1}(X)}$, we have $T^{j_{q+1}+m}y\in
\overline{T^{m+j_{q+1}+1}(X)}$ and therefore
$$
x+T^{j_{q+1}+m}y\in T^{m+j_q+1}(X)+\overline{T^{m+j_{q+1}+1}(X)}.
$$
Hence there exists $w_q\in X$ for which
$$
x+T^{j_{q+1}+m}y+T^{j_{q+1}+m+1}w_q\in
\overline{T^{m+j_{q+1}+1}(X)}.
$$
Then we can choose $v_q\in X$ such that
$$
\|x+T^{j_{q+1}+m}y+T^{j_{q+1}+m+1}w_q-T^{m+j_{q+1}+1}v_q\|\leq
2^{-q}.
$$
Denoting $u_{j_q}=w_q-T^{j_{q+1}-j_q}v_q$, we obtain
$$
x+T^{j_q+m}(y+Tu_{j_q})\in \overline{T^{m+j_{q+1}+1}(X)}\ \
\text{and}\ \ \|x+T^{j_q+m}(y+Tu_{j_q})\|\leq 2^{-q},
$$
which is exactly (\ref{indu}) for $n=q$. The construction of the
sequence $\{u_{j_n}\}_{n\in\zp}$, satisfying (\ref{indu}) for each
$n\in\zp$ is complete. From the inequality in (\ref{indu}) it
follows that the partial sums of the series
$\sum\limits_{n=0}^\infty T^{j_n+m}(y+Tu_{j_n})$ converge to zero
with respect to the norm of $X$. Hence $T^m(y+Tu_0)$ belongs to the
closed linear span of
$$
\{T^{m+n}(y+Tu_{j_n}):n\geq 1\}=\{T^{m+n}(y+Tu_n):n\in
A\setminus\{0\}\}
$$
as required. \square

\lemma{7.7}Let $T$ be an injective  bounded linear operator on a
Banach space $X$ such that $\overline{T^{m+1}(X)}$ has finite
codimension in $\overline{T^{m}(X)}$ for some $m\in\zp$. Then
$B(T,X)=A(T,X)$. \epr

\proof The inclusion $B(T,X)\subseteq A(T,X)$ follows from
Lemma~7.2. Let $y\in A(T,X)\setminus T(X)$. Since
$T^{-n}(\overline{T^{n+1}(X)})$ is an increasing sequence of linear
subspaces of $X$, there exists $q\in\zp$ such that $q\geq m$ and
$y\in T^{-q}(\overline{T^{q+1}(X)})$. According to Lemma~7.4
$\overline{T^{q+1}(X)}$ has finite codimension in
$\overline{T^{q}(X)}$.

Choose a strictly increasing sequence $\{m_k\}_{k\in\zp}$ of
positive integers and a sequence $\{A_k\}_{k\in\zp}$ of infinite
subsets of $\zp$ such that $m_0\geq q$, $0\in A_k$ for each
$k\in\zp$ and the sets $B_k=m_k+A_k=\{m_k+n:n\in A_k\}$ are
disjoint. For instance, we can choose a strictly increasing sequence
$\{m_n\}_{n\in\zp}$ of prime numbers such that $m_0\geq q$ and take
$A_k=\{m_k^l-m_k:l=1,2,\dots\}$. According to Lemma~7.6, for any
$k\in\zp$, there exists a sequence $\{u_n\}_{n\in B_k}$ of elements
of $X$ such that $T^q(y+Tu_{m_k})$ belongs to the closure of
$\li\{T^{q+j}(y+Tu_{m_k+j}):j\in A_k\setminus\{0\}\}$. For
$m\in\zp\setminus\bigcup\limits_{k=0}^\infty B_k$, we put $u_m=0$.
Since $T$ is continuous, we see that $T^{m_k-l}(y+Tu_{m_k})$ belongs
to the closed linear span of $\{T^{m_k-l+j}(y+Tu_{m_k+j}):j\geq 1\}$
if $q+l\leq m_k$. From Lemma~7.3 it follows now that the sequence
$\{\Phi_n=\Phi_n^{T,\kappa}\}_{n\in\zp}$ of linear functionals on
$\ssub XT(\kappa)$ for $\kappa=\{y+Tu_n\}_{n\in\zp}$ is not weakly
bounded. Hence $y\in B(T,X)$. \square

\lemma{7.8}Let $T$ and $S$ be injective bounded linear operators on
Banach spaces $X$ and $Y$ respectively and $G:X\to Y$ be an
invertible linear operator such that $SG=GT$. Suppose also that
$\overline{T^{m+1}(X)}$ has finite codimension in
$\overline{T^{m}(X)}$ for some $m\in\zp$. Then $G(A(T,X))=A(S,Y)$,
$A(T,X)$ is a closed finite codimensional linear subspace of $X$ and
$A(S,Y)$ is a closed linear subspace of $Y$. \epr

\proof Since the definition of the set $B(T,X)$ is "algebraic" we
have $G(B(T,X))=B(S,Y)$. By Lemma~7.7, we have $B(T,X)=A(T,X)$.
According to Lemma~7.4 $A(T,X)$ is a linear subspace of $X$ of
finite codimension, the increasing sequence
$T^{-n}(\overline{T^{n+1}(X)})$ of closed linear subspaces
stabilizes and $A(T,X)=T^{-n}(\overline{T^{n+1}(X)})$ for
sufficiently large $n$. In particular, $A(T,X)$ is closed in $X$.
Since $G(A(T,X))=G(B(T,X))=B(S,Y)$, we see that $B(S,Y)$ is a finite
codimensional linear subspace of $Y$. By Lemma~7.2, $A(S,Y)\supseteq
B(S,Y)$ and therefore $A(S,Y)$ has finite codimension in $Y$.
Therefore the increasing sequence $S^{-n}(\overline{S^{n+1}(Y)})$ of
closed linear subspaces of $Y$ stabilizes and eventually has finite
codimension. From Lemma~7.4 it follows that the codimension of
$\overline{S^{n+1}(Y)}$ in $\overline{S^{n}(Y)}$ is finite for
sufficiently large $n$. Applying Lemma~7.7 once again, we obtain
that $B(S,Y)=A(S,Y)$. Since by Lemma~7.4
$A(S,Y)=S^{-n}(\overline{S^{n+1}(Y)})$ for sufficiently large $n$,
we see that $A(S,Y)=B(S,Y)$ is a closed finite codimensional
subspace of $Y$. Finally, putting the above equalities together, we
see that $G(A(T,X))=G(B(T,X))=B(S,Y)=A(S,Y)$. \square

\subsection{Proof of Theorem~1.5}

Let $S$ be a bounded linear operator on a Banach space $Y$ and
$G:X\to Y$ be an invertible linear operator such that $SG=GT$.
According to Proposition~1.1 it suffices to prove that $G$ is
continuous. Injectivity of $T$ implies injectivity of $S$.

Consider $H_n\subseteq X$ defined inductively by the formulas
$$
H_0=X,\quad H_{k+1}=A(T,H_k)=\bigcup\limits_{n=0}^\infty
T^{-n}(\overline{T^{n+1}(H_k)})\quad \text{for $k\in\zp$.}
$$
According to Lemma~7.5 $H_k$ are closed finite codimensional linear
subspaces of $X$ and $\bigcap\limits_{k=0}^\infty H_k=\{0\}$.
Applying Lemma~7.8 consecutively to the restrictions of $T$ to the
invariant subspaces $H_k$, we obtain that $G(H_k)$ are closed in
$Y$. From Lemma~7.1 it now follows that $G$ is continuous. The proof
is complete.

\section{Concluding remarks}

{\bf1. }Let $\H=\ell_2\oplus L_2[0,1]$ and $T:\H\to \H$ be the
operator defined by the formula
$$
T(x\oplus f)=Ax\oplus Vf,
$$
where $V$ is the Volterra operator acting on $L_2[0,1]$ and
$A:\ell_2\to\ell_2$ be the weighted forward shift given by
$Ae_n=\frac{e_{n+1}}{n+1}$, $\{e_n\}_{n\in\zp}$ being the standard
orthonormal basis in $\ell_2$. Clearly $T$ is bounded injective and
quasinilpotent.

Using Proposition~4.4 and Lemma~4.6 it is easy to see that
$\ssub\dim\RR \ssub\H T=\ccc$ and $\ssub\dim\C \H/T(\H)=\ccc$. On
the other hand $\overline{T^{n+1}(X)}\neq \overline{T^{n}(X)}$ for
each $n\in\zp$.  According to Theorem~5.8, $T$ is not tame and
therefore not similar to the Volterra operator acting on $C[0,1]$.
It also does not determine the topology of $\H$ since it $T$ is the
direct sum of the Volterra operator $V$ acting on $L_2[0,1]$ with
another operator and $V$ does not determine the topology of
$L_2[0,1]$ according to Proposition~1.1  and Theorem~1.2. It worth
noting that $T$ satisfies (\ref{152}) and does not satisfy
(\ref{151}). This leads to the following natural question.

\smallskip

\noindent {\bf Problem 1.} \sl  Let $T$ be a bounded injective
linear operator on a Banach space $X$ such that $(\ref{151})$ is
satisfied. Is it true that $T$ determines the topology of $X$? \rm

\smallskip

\noindent{\bf 2. }A characterization of similarity of linear
operators on finite dimensional vector spaces is provided by the
Jordan block decomposition theorem, from which it follows that the
spectrum $\sigma(T)$ and the dimensions of ker$\,(T-\lambda I)^n$
for $\lambda\in\sigma(T)$, $n\in\zp$ are similarity invariants,
which determine $T$ up to similarity. This is obviously not true in
infinite dimensional case, when we have that the codimensions of
$(T-\lambda I)^n(X)$ and of $\bigcap\limits_{k=0}^\infty (T-\lambda
I)^k(X)$, which are also similarity invariants, are not determined
by the first family of invariants. The natural conjecture that these
dimensions and codimensions altogether form a set of invariants
determining $T$ up to similarity fails dismally as follows from
Theorem~2.2, as well as from the above example. This leads us to the
following problem.

\smallskip

\noindent {\bf Problem 2.} \sl Characterize similarity of linear
operators on infinite dimensional vector spaces. \rm

\smallskip

\noindent{\bf 3. }It worth noting that under the continuum
hypothesis Corollary~4.3, Proposition~4.4, Lemma~4.5 and Lemma~4.6
become trivial consequences of the Baire Theorem. However if one
does not assume the continuum hypothesis the Baire category argument
fails, since it is compatible with ZFC\footnote{ZFC stands for the
Zermelo--Frenkel axioms plus the axiom of choice, see, for instance
\cite{set}.} that a separable infinite dimensional Banach space is a
union of less than continuum of compact subsets.

\noindent{\bf 4. }The class of operators satisfying the conditions
(\ref{151}) and (\ref{152}) of Theorem~1.5 is rather rich. For
instance, it contains bounded weighted forward shifts on the spaces
$\ell_p$, $1\leq p\leq\infty$ and the Volterra operator acting on
the Hardy spaces $\H^p$ of the unit disk for $1\leq p\leq \infty$.
In particular, it contains plenty of injective quasinilpotent
operators. This class is also closed under finite powers and finite
direct sums. It worth noting that (\ref{151}) and (\ref{152}) imply
that there are countably many continuous linear functionals on $X$
separating points and therefore the algebraic dimension of $X$ is
$\ccc$.

\noindent{\bf 5. }The following theorem is proved in \cite{ja1}.

\theorem{J} Let $A$ be a commutative unital semisimple Banach
algebra and $a\in A$. Then the multiplication operator $M_a:A\to A$,
$M_a x=ax$ determines the topology of $A$ if and only if for any
$\lambda\in\C$ either $a-\lambda$ is not a zero divisor or
$a-\lambda$ is a zero divisor and the codimension of $(M_a-\lambda
I)(A)$ in $A$ is finite. \epr

The classes of operators on Banach spaces determining the topology,
provided by Theorem~J and by Theorem~1.5 do not cover each other
(although they do intersect). For instance, the multiplication
operator $T:C[0,1]\to C[0,1]$, $Tf(x)=xf(x)$ determines the topology
of $C[0,1]$ according to Theorem~J and it does not satisfy
(\ref{151}). On the other hand a multiplication operator $M_a$ with
$a\neq 0$ on a commutative unital semisimple Banach algebra is never
quasinilpotent, while Theorem~1.5 can be applied to certain
quasinilpotent operators like the Volterra operator acting on the
Hardy space $\H^p$ of the unit disk. It would be interesting to find
a unified approach, generalizing Theorems~J and~1.5 simultaneously.

\bigskip

{\bf Acknowledgements.} \ Partially supported by Plan Nacional I+D+I
Grant BFM2003-00034, Junta de Andaluc\'{\i}a FQM-260 and British
Engineering and Physical Research Council Grant GR/T25552/01.


\end{document}